\documentclass[11pt,reqno]{article}

\usepackage{amsfonts}
\usepackage{amsmath}
\usepackage{amsthm}
\usepackage{amssymb}
\usepackage{mathrsfs}
\usepackage{amstext}
\usepackage{graphicx}
\font\msbm=msbm10

\numberwithin{equation}{section}

\theoremstyle{plain}
\newtheorem{Theorem}{Theorem}[section]
\newtheorem{lemma}[Theorem]{Lemma}
\newtheorem{corollary}[Theorem]{Corollary}

\newtheorem{proposition}[Theorem]{Proposition}
\newtheorem{definition}[Theorem]{Definition}

\newtheorem{remark}[Theorem]{Remark}
\def\mathbb#1{\hbox{\msbm{#1}}}

\newcommand{\re}{\ensuremath{\mathbb{R}}}\newcommand{\N}{\ensuremath{\mathbb{N}}}
\newcommand{\zz}{\ensuremath{\mathbb{Z}}}

\newcommand{\com}{\ensuremath{\mathbb{C}}}
\newcommand{\Z}{{\ensuremath{\zz}^d}}
\newcommand{\n}{\ensuremath{{\N}_0}}
\newcommand{\R}{\ensuremath{{\re}^d}}

\newcommand{\bproof}{\noindent {\bf Proof.}\, }
\newcommand{\eproof}{\hspace*{\fill} \rule{3mm}{3mm}}

\newcommand{\CoY}{\ensuremath{\mbox{Co}Y}}

\newcommand{\supp}{{\rm supp \, }}

\newcommand{\frei}{\hspace*{\fill}}

\newcommand{\cd}{\ensuremath{\mathcal D}}

\newcommand{\cx}{\ensuremath{\mathcal X}}

\newcommand{\cn}{\ensuremath{\mathcal N}}
\newcommand{\cf}{\ensuremath{\mathcal F}}
\newcommand{\cfi}{\ensuremath{{\cf}^{-1}}}

\newcommand{\cg}{\ensuremath{\mathcal G}}

\newlength{\fixboxwidth}
\newcommand{\be}{\begin{equation}}
\newcommand{\ee}{\end{equation}}
\newcommand{\beq}{\begin{eqnarray}}
\newcommand{\beqq}{\begin{eqnarray*}}
\newcommand{\eeq}{\end{eqnarray}}
\newcommand{\eeqq}{\end{eqnarray*}}

\newcommand{\dotspace}[2]{\ensuremath{\dot{#1}^{#2}_{p,q}}(\mathcal{G})}

\newcommand{\Dil}{\ensuremath{\mathcal{D}}}

\begin{document}
\title{Continuous characterizations of Besov-Lizorkin-Triebel spaces
and new interpretations as coorbits}

\author{Tino Ullrich}

\maketitle

\begin{abstract} We give characterizations for homogeneous
and inhomogeneous Besov-Lizorkin-Triebel spaces \cite{Tr83,Tr92,Tr06}
in terms of continuous local means for the full range of parameters. In
particular, we prove characterizations using tent spaces (Lusin functions) and
spaces involving the Peetre maximal function in order to apply the classical
coorbit space theory due to Feichtinger and Gr\"ochenig 
\cite{FeGr86, FeGr89a, FeGr89b, Gr88, Gr91}. This
results in atomic decompositions and wavelet bases for homogeneous
spaces. In particular we give sufficient conditions for suitable
wavelets in terms of moment, decay and smoothness conditions.
\end{abstract}

\begin{tabbing}
Peetre maximal
function, tent spaces, atomic decompositions,
wavelet bases.
\end{tabbing}

\section{Introduction}
This paper deals with Besov-Lizorkin-Triebel spaces $\dot{B}^s_{p,q}(\R)$ and $\dot{F}^s_{p,q}(\R)$
on the Euclidean space $\R$ and their interpretation as coorbits. 
For this purpose  
we prove a number of characterizations for homogeneous and inhomogeneous spaces 
for the full range of parameters. Classically introduced in Triebel's monograph \cite[2.3.1]{Tr83} by means of a dyadic decomposition
of unity, we use more general building blocks
and provide in addition continuous characterizations in terms of Lusin and maximal functions. Equivalent 
(quasi-)normings of this kind were first given by Triebel in \cite{Tr88}. His proofs use in an essential way the fact that the
function under consideration belongs to the respective space. Therefore, 
the obtained equivalent (quasi-) norms could not yet be considered as a definition or characterization of the space.
Later on, Triebel was able to solve this problem partly in his monograph \cite[2.4.2, 2.5.1]{Tr92} by restricting to the Banach space
case. Afterwards, Rychkov \cite{Ry99a} completed the picture 
by simplifying a method due to Bui, Paluszy\'nski, and Taibleson \cite{BuPaTa96,BuPaTa97}.
However, \cite{Ry99a} contains some problematic arguments. One aim of the 
present paper is to provide a complete and self-contained reference for general 
characterizations of discrete and continuous type by avoiding these arguments. 
We use a variant of a method from Rychkov's subsequent papers \cite{Ry99b,Ry01} which is 
originally due to Str{\"o}mberg and Torchinsky developed in their monograph \cite[Chapt.
5]{StTo89}.

In a different language the results can be interpreted in terms of
the continuous wavelet transform (see Appendix \ref{SectCWT}) belonging to a function space on
the $ax+b$-group $\cg$. Spaces on $\cg$ considered here are
mixed norm spaces like tent spaces \cite{CoMeSt85} as well as Peetre type spaces. The latter are indeed new
and received their name from the fact that quantities related to the classical Peetre
maximal function are involved. This leads to the main intention of the paper.
We use the established characterizations for the homogeneous spaces in
order to embed them in the abstract
framework of coorbit space theory originally due to Feichtinger
and Gr\"ochenig \cite{FeGr86, FeGr89a, FeGr89b, Gr88, Gr91} in the 80s. This connection was already
observed by them in \cite{FeGr86,Gr88,Gr91}. They worked with Triebel's equivalent
continuous normings from \cite{Tr88} and the results on tent spaces which were introduced more or less
at the same time by Coifman, Meyer, Stein \cite{CoMeSt85} to
interpret Lizorkin-Triebel spaces as coorbits. On the one hand the present paper gives a late
justification and on the other hand we observe that Peetre type spaces on $\cg$ are a
much better choice  for this issue. Their two-sided translation invariance is
immediate and much more transparent as we will show in Section \ref{PeFS}. Furthermore,
generalizations in different directions are now possible. In a
forthcoming paper we will show how to apply a generalized coorbit
space theory due to Fornasier and Rauhut \cite{fora05} in order to
recover inhomogeneous spaces based on the characterizations given here. 
Moreover, the extension of the results
to quasi-Banach spaces using a
theory developed by Rauhut in \cite{ra05-3,ra05-4} is possible.

Once we have interpreted classical homogeneous Besov-Lizorkin-Triebel spaces as certain
coorbits, we are able to benefit from the achievements of the
abstract theory in \cite{FeGr86, FeGr89a, FeGr89b, Gr88, Gr91}. 
The main feature is a powerful
discretization machinery which leads in an abstract universal way to
atomic decompositions. We are now able to apply this method which results in 
atomic decompositions and wavelet bases for homogeneous spaces. More precisely, sufficient conditions
in terms of vanishing moments, decay, and smoothness properties of
the respective wavelet function are given. Compact support of the
used atoms does not play any role here. In particular, we specify
the order of a suitable orthonormal spline wavelet system depending
on the parameters of the respective space.

The paper is organized as follows. After giving some preliminaries
we start in Section 2 with the definition of classical
Besov-Lizorkin-Triebel spaces and their characterization via
continuous local means. In Section 3 we give a brief introduction to
abstract coorbit space theory which is applied in Section 4 on the
$ax+b$-group $\mathcal{G}$. We recover the homogeneous spaces from Section 2 as
coorbits of certain spaces on $\mathcal{G}$. Finally, several
discretization results in terms of atomic decompositions and wavelet isomorphisms 
are established. The underlying decay result of the continuous wavelet transform 
and some basic facts about orthonormal wavelet bases are shifted to the appendix.

{\bf Acknowledgement:} The author would like to thank Holger Rauhut, Martin Sch\"afer, Benjamin Scharf, and Hans Triebel 
for valuable discussions, a critical reading of preliminary versions of this manuscript and for several hints
how to improve it.

\subsection{Notation}
Let us first introduce some basic notation. The symbols
$\re, \com, \N, \n$ and $\zz$ denote the real numbers, complex
numbers, natural numbers, natural numbers including 0 and the
integers. The dimension of the underlying Euclidean space for
function spaces is denoted by $d$, its elements will be denoted by
$x,y,z,...$ and $|x|$ is used for the Euclidean norm. We will use
$|k|_1$ for the $\ell_1^d$-norm of a vector $k$. For a
multi-index $\bar{\alpha}$ and $x\in \R$ we write
$$
    x^{\bar{\alpha}} = x_1^{\alpha_1}\cdot...\cdot x_{d}^{\alpha_d}
$$
and define the differential operators $D^{\bar{\alpha}}$ and $\Delta$ by
$$
    D^{\bar{\alpha}} = \frac{\partial^{|\bar{\alpha}|_1}}{\partial x_1^{\alpha_1}\cdots
    \partial x_d^{\alpha_d}}\quad\mbox{and}\quad
    \Delta = \sum\limits_{k=1}^d \frac{\partial^2}{\partial x_k^2}\,.
$$
If $X$ is a (quasi-)Banach space and $f\in X$ we use $\|f|X\|$ or
simply $\|f\|$ for its\\ (quasi-)norm. The space of linear continuous mappings from 
$X$ to $Y$ is denoted by $\mathcal{L}(X,Y)$ or simply $\mathcal{L}(X)$ if $X=Y$. Operator 
(quasi-)norms of $A \in \mathcal{L}(X,Y)$ are denoted by $\|A:X\to Y\|$, or simply by
$\|A\|$. As usual, the letter $c$ denotes a constant, which may vary
from line to line but is always independent of $f$, unless the
opposite is explicitly stated. We also use the notation
$a\lesssim b$ if there exists a constant $c>0$ (independent of the
context dependent relevant parameters) such that $a \le c\,b$. If
$a\lesssim b$ and $b \lesssim a$ we will write $a \asymp b$\,.

\section{Function spaces on $\R$}
\label{clFS}
\subsection{Vector valued Lebesgue spaces}
The space $L_p(\R)$, $0<p\leq \infty$, denotes the collection of 
complex-valued functions (equivalence classes) with finite (quasi-)norm
$$
    \|f|L_p(\R)\| = \Big(\,\int\limits_{\R}|f(x)|^p\,dx\Big)^{1/p}\,,
$$ with the usual modification if $p=\infty$. The Hilbert space $L_2(\R)$
plays a separate role for our purpose (Section \ref{Classcoo}). Having a sequence of complex-valued functions
$\{f_{k}\}_{k \in I}$ on $\R$, where $I$ is a countable index set,
we put
$$
   \|\{f_{k}\}_k|\ell_q(L_p(\R))\| = \Big(\sum\limits_{k\in I}\|f_{k}|L_p(\R)\|^q\Big)^{1/q}
$$
and
$$
  \|\{f_{k}\}_k|L_p(\ell_q,\R)\| = \Big\|\Big(\sum\limits_{k\in I}|f_{k}(x)|^q\Big)^{1/q}\Big|L_p(\R)\Big\|\,,
$$
where we modify appropriately in the case $q=\infty$.\\

\subsection{Maximal functions}
\label{sectmax}
For a locally integrable function $f$ we denote by $Mf(x)$ the
Hardy-Littlewood maximal function defined by
\begin{equation}\nonumber
  (Mf)(x) = \sup\limits_{x\in Q} \frac{1}{|Q|}\int\limits_{Q}\,|f(y)|\,dy\quad,\quad x\in\R
  \,,
\end{equation}
where the supremum is taken over all cubes centered at $x$ with
sides parallel to the coordinate axes. The following theorem is due
to Fefferman and Stein \cite{FeSt71}.
\begin{Theorem}\label{feffstein}\rm For $1<p<\infty$ and $1 <q \leq \infty$ there exists a constant $c>0$, such that
  $$
      \|\{Mf_{k}\}_k|L_p(\ell_q,\R)\| \leq c\|\{f_{k}\}_k|L_p(\ell_q,\R)\|
  $$
  holds for all sequences $\{f_k\}_{k\in \zz}$ of locally Lebesgue-integrable functions on $\R$.
\end{Theorem}

Let us recall the classical Peetre maximal function, introduced in
\cite{Pe75}\,. Given a sequence of functions $\{\Psi_k\}_{k\in \N}
\subset \mathcal{S}(\R)$, a tempered distribution $f\in \mathcal{S}'(\R)$ and a positive
number $a>0$ we define the system of maximal functions
\begin{equation}\nonumber
    (\Psi_{k}^{\ast}f)_{a}(x) = \sup\limits_{y \in \R}\frac{|(\Psi_{k}\ast f)(x+y)|}
    {(1+2^k|y|)^{a}}\quad,\quad x\in \R, k\in \zz\,.
\end{equation}
Since $(\Psi_{k}\ast f)(y)$ makes sense pointwise (see the following
paragraph) everything is well-defined. However, the value
``$\infty$'' is also possible for $(\Psi_k^\ast f)_a(x)$. This
was the reason for the problematic arguments in \cite{Ry99a} mentioned in the introduction.
We will often use
dilates $\Psi_k(x) = 2^{kd}\Psi(2^kx)$ of a fixed function $\Psi \in
\mathcal{S}(\R)$, where $\Psi_0(x)$ might be given by a separate function. Also
continuous dilates are needed. Let the operator
$\cd_t^{L_p}$, $t>0$, generate the $p$-normalized dilates of a function
$\Psi$ given by $\cd_t^{L_p}\Psi := t^{-d/p}\Psi(t^{-1}\cdot)$. If
$p=1$ we omit the super index and use additionally $\Psi_t:=\cd_t
\Psi:=\cd_t^{L_1}\Psi$. We define $(\Psi^{\ast}_t f)_a(x)$ by
\begin{equation}\label{Peetrefunct}
    (\Psi^{\ast}_t f)_{a}(x) = \sup\limits_{y \in \R}\frac{|(\Psi_t \ast f)(x+y)|}
    {(1+|y|/t)^{a}}\quad,\quad x\in \R\,,t>0\,.
\end{equation}
We will refer to this construction later on. It turned out that this
maximal function construction can be used to interpret classical
smoothness spaces as coorbits of certain Banach function spaces on the
$ax+b$-group, see Section \ref{ax+b}.

\subsection{Tempered distributions, Fourier transform}

As usual $\mathcal{S}(\re^d)$ is used for
the locally convex space of rapidly decreasing infinitely
differentiable functions on $\re^d$ where its topology is generated
by the family of semi-norms
$$
    \|\varphi\|_{k,\ell} = \sup\limits_{x\in \re^d,|\bar{\alpha}|_1 \leq \ell}
    |D^{\bar{\alpha}}\varphi(x)|(1+|x|)^k\quad,\quad \varphi\in \mathcal{S}(\re^d)\,,~k,\ell\in \n\,.
$$
The space $\mathcal{S}'(\re^d)$, the topological dual of $\mathcal{S}(\R)$, is also referred as the set of tempered
distributions on $\re^d$. Indeed, a linear mapping $f:\mathcal{S}(\re^d) \to \com$ belongs
to $\mathcal{S}'(\re^d)$ if and only if there exist numbers $k,\ell \in \n$
and a constant $c = c_f$ such that
\begin{equation}\label{eq100}
    |f(\varphi)| \leq c_f\sup\limits_{x\in \R,|\bar{\alpha}|_1 \leq \ell}
    |D^{\bar{\alpha}}\varphi(x)|(1+|x|)^k
\end{equation}
for all $\varphi\in \mathcal{S}(\re^d)$. The space $\mathcal{S}'(\re^d)$ is
equipped with the weak$^{\ast}$-topology.

The convolution $\varphi\ast \psi$ of two
integrable (square integrable) functions $\varphi, \psi$ is defined via the integral
\begin{equation}\label{conv}
    (\varphi \ast \psi)(x) = \int\limits_{\R} \varphi(x-y)\psi(y)\,dy\,.
\end{equation}
If $\varphi,\psi \in \mathcal{S}(\R)$ then \eqref{conv} still belongs to
$\mathcal{S}(\R)$. The convolution can be extended to $\mathcal{S}(\R)\times \mathcal{S}'(\R)$ via
$(\varphi\ast f)(x) = f(\varphi(x-\cdot))$. It makes sense pointwise and is 
a $C^{\infty}$-function in $\R$ of at most polynomial growth. 

As usual the Fourier transform defined on both $\mathcal{S}(\R)$ and $\mathcal{S}'(\R$)
is given by $(\cf f)(\varphi) := f(\cf \varphi)$, where
$\,f\in \mathcal{S}'(\R), \varphi \in \mathcal{S}(\R)$, and
$$
\cf \varphi(\xi) := (2\pi)^{-d/2}\int\limits_{\R} e^{-ix\cdot
\xi}\varphi(x)\,dx.
$$
The mapping $\cf$ is a bijection (in both cases) and its inverse is
given by $\cf^{-1}\varphi = \cf\varphi(-\cdot)$.

In order to deal with homogeneous spaces we need to define the subset $\mathcal{S}_0(\R)\subset \mathcal{S}(\R)$.
Following \cite[Chapt. 5]{Tr83} we put 
$$
    \mathcal{S}_0(\R) = \{\varphi\in \mathcal{S}(\R)~:~D^{\bar{\alpha}}(\cf \varphi)(0)=0\quad
    \mbox{for every multi-index }\bar{\alpha} \in \n^d\}\,.
$$
The set $\mathcal{S}_0'(\R)$ denotes the topological dual of $\mathcal{S}_0(\R)$\,.
If $f\in \mathcal{S}'(\R)$, the restriction of $f$ to $\mathcal{S}_0(\R)$ clearly belongs to $\mathcal{S}_0'(\R)$\,.
Furthermore, if $P(x)$ is an arbitrary
polynomial in $\R$, we have $(f+P(\cdot))(\varphi) = f(\varphi)$
for every $\varphi\in \mathcal{S}_0(\R)$. Conversely, if $f\in \mathcal{S}_0'(\R)$, then $f$ can be extended
from $\mathcal{S}_0(\R)$ to $\mathcal{S}(\R)$, i.e., to an element of $\mathcal{S}'(\R)$\,. However, this fact is not trivial
and makes use of the Hahn-Banach theorem in locally convex topological vector spaces.
We may identify $\mathcal{S}_0'(\R)$ with the factor space $\mathcal{S}'(\R)/\mathcal{P}(\R)$, since two
different extensions differ by a polynomial\,.

\subsection{Besov-Lizorkin-Triebel spaces}

Let us first introduce the concept of a dyadic decomposition of
unity, see also \cite[2.3.1]{Tr83}.

\begin{definition} {\em(a)} Let $\Phi(\R)$ be the collection of all systems $\{\varphi_j(x)\}_{j\in
\n} \subset \mathcal{S}(\R)$ with the following properties
\begin{description}
    \item(i) $\varphi_j(x) = \varphi(2^{-j}x)\quad,\quad j\in
    \N$\quad,
    \item(ii) $\supp \varphi_0 \subset \{x\in \R~:~|x|\leq
    2\}\quad,\quad
    \supp \varphi \subset \{x\in \R~:~1/2 \leq |x| \leq 2\}\quad,\quad$and
    \item(iii) $\sum\limits_{j=0}^{\infty} \varphi_j(x) = 1$ for every $x\in \R$\,.
\end{description}
{\em(b)} Moreover, $\dot{\Phi}(\R)$ denotes the collection of all
systems $\{\varphi_j(x)\}_{j\in \zz} \subset \mathcal{S}(\R)$ with the
following properties
\begin{description}
    \item(i) $\varphi_j(x) = \varphi(2^{-j}x)\quad,\quad j\in
    \zz\quad,$
    \item(ii) $\supp \varphi = \{x\in \R~:~1/2 \leq |x| \leq 2\}\quad,\quad$and
    \item(iii) $\sum\limits_{j=-\infty}^{\infty}\varphi_j = 1$ for every $x\in \R
    \setminus\{0\}$\,.
\end{description}
\end{definition}

\begin{remark}\label{remunity} If we take $\varphi_0 \in \mathcal{S}(\R)$ satisfying
$$
     \varphi_0(x) = \left\{\begin{array}{lcr}
     1&:&|x| \leq 1\\
     0&:&|x|>2
    \end{array}\right.
$$
and define $\varphi(x) = \varphi_0(x)-\varphi_0(2x)$ then the system
$\{\varphi_j(x)\}_{j\in \n}$ belongs to $\Phi(\R)$ and the system $\{\varphi_j(x)\}_{j\in \zz}$
with $\varphi_0 := \varphi$ belongs to $\dot{\Phi}(\R)$.
\end{remark}

Now we are ready for the definition of the Besov and
Lizorkin-Triebel spaces. See for instance \cite[2.3.1]{Tr83} for details and further properties.

\begin{definition}\label{inhom} Let $\{\varphi_j(x)\}_{j=0}^{\infty} \in \Phi(\R)$ and
$\Phi_j = \cf^{-1}\varphi_j$, $j\in \n$.
Let further $-\infty<s<\infty$ and $0<q\leq \infty$.
\begin{description}
 \item(i) If $0<p\leq \infty$ then
 $$
    B^s_{p,q}(\R) = \Big\{f\in \mathcal{S}'(\R):
    \|f|B^s_{p,q}(\R)\| = \Big(\sum\limits_{j=0}^{\infty}2^{jsq}
    \|\Phi_j\ast f|L_p(\R)\|^q\Big)^{1/q}<\infty\Big\}\,.
 $$
 \item(ii) If $0<p<\infty$ then
 $$
    F^s_{p,q}(\R) = \Big\{f\in \mathcal{S}'(\R):
    \|f|F^s_{p,q}(\R)\| = \Big\|\Big(\sum\limits_{j=0}^{\infty}2^{jsq}
    |(\Phi_j \ast f)(x)|^q\Big)^{1/q}|L_p(\R)\Big\|<\infty\Big\}\,.
 $$
\end{description}
In case $q=\infty$ we replace the sum by a supremum in both cases.
\end{definition}
\noindent The homogeneous counterparts are defined as follows. For details, further properties and how to deal with 
occurring technicalities we refer to \cite[Chapt. 5]{Tr83}.

\begin{definition} Let $\{\varphi_j(x)\}_{j \in \zz} \in \dot{\Phi}(\R)$ and
$\Phi_j = \cf^{-1}\varphi_j$.
Let further $-\infty<s<\infty$ and $0<q\leq \infty$.
\begin{description}
 \item(i) If $0<p\leq \infty$ then
 $$
    \dot{B}^s_{p,q}(\R) = \Big\{f\in \mathcal{S}_0'(\R):
    \|f|\dot{B}^s_{p,q}(\R)\| = \Big(\sum\limits_{j=-\infty}^{\infty}2^{jsq}
    \|\Phi_j\ast f|L_p(\R)\|^q\Big)^{1/q}<\infty\Big\}\,.
 $$
 \item(ii) If $0<p<\infty$ then
 $$
    \dot{F}^s_{p,q}(\R) = \Big\{f\in \mathcal{S}_0'(\R):
    \|f|\dot{F}^s_{p,q}(\R)\| = \Big\|\Big(\sum\limits_{j=-\infty}^{\infty}2^{jsq}
    |(\Phi_j \ast f)(x)|^q\Big)^{1/q}|L_p(\R)\Big\|<\infty\Big\}\,.
 $$
\end{description}
In case $q=\infty$ we replace the sum by a supremum in both cases.
\end{definition}

\subsection{Inhomogeneous spaces}
Essential for the sequel are functions $\Phi_0, \Phi \in \mathcal{S}(\R)$ satisfying
\begin{equation}\label{condphi1}
\begin{split}
  &|(\cf \Phi_0)(x)| > 0 \quad \mbox{ on }\quad \{|x| < 2\varepsilon\}\\
  &|(\cf \Phi)(x)| > 0 \quad\mbox{ on }\quad \{\varepsilon/2<|x|< 2\varepsilon\}\,,
\end{split}
\end{equation}
for some $\varepsilon>0$, and
\begin{equation}\label{condphi2}
   D^{\bar{\alpha}} (\cf \Phi)(0) = 0\quad\mbox{for all}\quad
   |\bar{\alpha}|_1 \leq R.
\end{equation}
We will call the functions $\Phi_0$ and $\Phi$ kernels for local
means. Recall that $\Phi_k = 2^{kd}\Phi(2^k\cdot)$, $k\in \N$, and
$\Psi_t = \Dil_t \Psi$. The upcoming four theorems represent the main results of the first part
of the paper. 

\begin{Theorem}\label{main1} Let $s\in \re$, $0<p<\infty$, $0<q\leq
\infty$, $a>d/\min\{p,q\}$ and $R+1> s$. Let further $\Phi_0,\Phi\in \mathcal{S}(\R)$ 
be given by \eqref{condphi1} and \eqref{condphi2}. Then the space $F^s_{p,q}(\R)$
can be characterized by
$$
    F^s_{p,q}(\R) = \{f\in \mathcal{S}'(\R)~:~\|f|F^s_{p,q}(\R)\|_i <
    \infty\}\quad,\quad i=1,...,5,
$$
where
\begin{eqnarray}
    \|f|F^s_{p,q}\|_1 &=& \|\Phi_0 \ast f|L_p(\R)
    \| + \Big\|\Big(\int\limits_{0}^1 t^{-sq}
    |(\Phi_t \ast f)(x)|^q
    \frac{dt}{t}\Big)^{1/q}|L_p(\R)\Big\|~,\label{CL}\\
    \|f|F^s_{p,q}\|_2 &=& \|(\Phi_0^{\ast}f)_a|L_p(\R)
    \| + \Big\|\Big(\int\limits_{0}^1 t^{-sq}
    \Big[\sup\limits_{z\in \R}\frac{|(\Phi_t \ast
    f)(x+z)|}{(1+|z|/t)^a}\Big]^q
    \frac{dt}{t}\Big)^{1/q}|L_p(\R)\Big\|~,\label{CP}\\
    \|f|F^s_{p,q}\|_3 &=& \|\Phi_0 \ast f|L_p(\R)\|
    + \Big\|\Big(\int\limits_{0}^1 t^{-sq} \int\limits_{|z| < t}
    |(\Phi_t \ast f)(x+z)|^q
    \frac{dt}{t^{d+1}}\Big)^{1/q}|L_p(\R)\Big\|~,\label{CT}\\
    \|f|F^s_{p,q}\|_4 &=& \Big\|\Big(\sum\limits_{k=0}^{\infty}2^{skq}
    \Big[\sup\limits_{z\in \R}\frac{|(\Phi_k \ast f)(x+z)|}{(1+2^k|z|)^a}\Big]^q\Big)^{1/q}|L_p(\R)\Big\|~,\label{DP}\\
    \|f|F^s_{p,q}\|_5 &=&
    \Big\|\Big(\sum\limits_{k=0}^{\infty}2^{skq}
    |(\Phi_k \ast f)(x)|^q\Big)^{1/q}|L_p(\R)\Big\|~\label{DL}
\end{eqnarray}
with the usual modification in case $q=\infty$. Furthermore, all quantities
$\|f|F^s_{p,q}(\R)\|_i$, $i=1,...,5$, are equivalent (quasi-)norms
in $F^s_{p,q}(\R)$\,.
\end{Theorem}
For the inhomogeneous Besov spaces we obtain the following. 
\begin{Theorem}Let $s\in \re$, $0<p,q\leq\infty$, $a>d/p$ and $R+1> s$. Let further $\Phi_0,\Phi\in \mathcal{S}(\R)$ 
be given by \eqref{condphi1} and \eqref{condphi2}. Then the space
$B^s_{p,q}(\R)$ can be characterized by
$$
    B^s_{p,q}(\R) = \{f\in \mathcal{S}'(\R)~:~\|f|B^s_{p,q}(\R)\|_i <
    \infty\}\quad,\quad i=1,...,4,
$$
where
\begin{eqnarray}
    \|f|B^s_{p,q}\|_1 &=& \|\Phi_0 \ast f|L_p(\R)
    \| + \Big(\int\limits_{0}^1 t^{-sq}
    \|(\Phi_t \ast f)(x)|L_p(\R)\|^q
    \frac{dt}{t}\Big)^{1/q}~,\nonumber\\
    \|f|B^s_{p,q}\|_2 &=& \|(\Phi_0^{\ast}f)_a|L_p(\R)
    \| + \Big(\int\limits_{0}^1 t^{-sq}
    \Big\|\sup\limits_{z\in \R}\frac{|(\Phi_t \ast
    f)(x+z)|}{(1+|z|/t)^a}|L_p(\R)\Big\|^q\frac{dt}{t}\Big)^{1/q},\nonumber\\
    \|f|B^s_{p,q}\|_3 &=&
    \Big(\sum\limits_{k=0}^{\infty}2^{skq}
    \Big\|\sup\limits_{z\in \R}\frac{|(\Phi_k \ast f)(x+z)|}{(1+2^k|z|)^a}
    |L_p(\R)\Big\|^q\Big)^{1/q}~,\nonumber\\
    \|f|B^s_{p,q}\|_4 &=&
    \Big(\sum\limits_{k=0}^{\infty}2^{skq}
    \|(\Phi_k \ast f)(x)|L_p(\R)\|^q\Big)^{1/q}~\nonumber
\end{eqnarray}
with the usual modification if $q=\infty$. Furthermore, all quantities
$\|f|B^s_{p,q}(\R)\|_i$, $i=1,...,4$, are equivalent quasi-norms in
$B^s_{p,q}(\R)$\,.
\end{Theorem}

\subsection{Homogeneous spaces}
The homogeneous spaces can be characterized similar. Here we do not
have a separate function $\Phi_0$ anymore. We put $\Phi_0 = \Phi$\,.

\begin{Theorem}\label{main1b} Let $s\in \re$, $0<p<\infty$, $0<q\leq
\infty$, $a>d/\min\{p,q\}$ and $R+1> s$. Let further $\Phi \in \mathcal{S}(\R)$ be
given by \eqref{condphi1} and \eqref{condphi2}. Then the space $\dot{F}^s_{p,q}(\R)$ can
be characterized by
$$
    \dot{F}^s_{p,q}(\R) = \{f\in \mathcal{S}_0'(\R)~:~\|f|\dot{F}^s_{p,q}(\R)\|_i <
    \infty\}\quad,\quad i=1,...,5,
$$
where\vspace{-0.3cm}
\begin{equation}\nonumber
  \begin{split}
    \|f|\dot{F}^s_{p,q}\|_1 &= \Big\|\Big(\int\limits_{0}^{\infty} t^{-sq}
    |(\Phi_t \ast f)(x)|^q
    \frac{dt}{t}\Big)^{1/q}|L_p(\R)\Big\|~,\\
    \|f|\dot{F}^s_{p,q}\|_2 &= \Big\|\Big(\int\limits_{0}^\infty t^{-sq}
    \Big[\sup\limits_{z\in \R}\frac{|(\Phi_t \ast
    f)(x+z)|}{(1+|z|/t)^a}\Big]^q
    \frac{dt}{t}\Big)^{1/q}|L_p(\R)\Big\|~,\\
    \|f|\dot{F}^s_{p,q}\|_3 &= \Big\|\Big(\int\limits_{0}^{\infty} t^{-sq} \int\limits_{|z| < t}
    |(\Phi_t \ast f)(x+z)|^q
    \frac{dt}{t^{d+1}}\Big)^{1/q}|L_p(\R)\Big\|~,\\
    \|f|\dot{F}^s_{p,q}\|_4 &=
    \Big\|\Big(\sum\limits_{k=-\infty}^{\infty}2^{skq}
    \Big[\sup\limits_{z\in \R}\frac{|(\Phi_k \ast f)(x+z)|}{(1+2^k|z|)^a}\Big]^q
    \Big)^{1/q}|L_p(\R)\Big\|~,\\
    \|f|\dot{F}^s_{p,q}\|_5 &=
    \Big\|\Big(\sum\limits_{k=-\infty}^{\infty}2^{skq}
    |(\Phi_k \ast f)(x)|^q\Big)^{1/q}|L_p(\R)\Big\|~
  \end{split}
\end{equation}
with the usual modification if $q=\infty$. Furthermore, all quantities
$\|f|\dot{F}^s_{p,q}(\R)\|_i$, $i=1,...,5$, are equivalent
quasi-norms in $\dot{F}^s_{p,q}(\R)$\,.
\end{Theorem}
For the homogeneous Besov spaces we obtain the following.
\begin{Theorem}\label{main1bb}Let $s\in \re$, $0<p,q\leq\infty$,
$a>d/p$ and $R+1> s$. Let further $\Phi \in \mathcal{S}(\R)$ be
given by \eqref{condphi1} and \eqref{condphi2}. Then the space
$\dot{B}^s_{p,q}(\R)$ can be characterized by
$$
    \dot{B}^s_{p,q}(\R) = \{f\in \mathcal{S}_0'(\R)~:~\|f|\dot{B}^s_{p,q}(\R)\|_i <
    \infty\}\quad,\quad i=1,...,4,
$$
where
\begin{eqnarray}
    \|f|\dot{B}^s_{p,q}\|_1 &=& \Big(\int\limits_{0}^\infty t^{-sq}
    \|(\Phi_t \ast f)(x)|L_p(\R)\|^q
    \frac{dt}{t}\Big)^{1/q}~,\nonumber\\
    \|f|\dot{B}^s_{p,q}\|_2 &=& \Big(\int\limits_{0}^\infty t^{-sq}
    \Big\|\sup\limits_{z\in \R}\frac{|(\Phi_t \ast
    f)(x+z)|}{(1+|z|/t)^a}|L_p(\R)\Big\|^q\frac{dt}{t}\Big)^{1/q}~\nonumber,\\
    \|f|\dot{B}^s_{p,q}\|_3 &=&
    \Big(\sum\limits_{k=-\infty}^{\infty}2^{skq}
    \Big\|\sup\limits_{z\in \R}\frac{|(\Phi_k \ast f)(x+z)|}{(1+2^k|z|)^a}|L_p(\R)
    \Big\|^q\Big)^{1/q}~,\nonumber\\
    \|f|\dot{B}^s_{p,q}\|_4 &=&
    \Big(\sum\limits_{k=-\infty}^{\infty}2^{skq}
    \|(\Phi_k \ast f)(x)|L_p(\R)\|^q\Big)^{1/q}~\nonumber
\end{eqnarray}
with the usual modification if $q=\infty$. Furthermore, all quantities
$\|f|\dot{B}^s_{p,q}(\R)\|_i$, $i=1,...,4$, are equivalent quasi-norms in
$\dot{B}^s_{p,q}(\R)$\,.
\end{Theorem}

\begin{remark} Observe, that the (quasi-)norms $\|\cdot|\dot{F}^s_{p,q}(\R)\|_3$ and $\|\cdot|F^s_{p,q}(\R)\|_3$ are characterizations
via Lusin functions, see 
\cite[2.4.5]{Tr92} and \cite[2.12.1]{Tr83} and the references given there. We will return to it later when defining tent spaces, see Definition \ref{sponG}
and \eqref{tent}.
\end{remark}

\subsection{Particular kernels}
For more details concerning particular choices for the kernels $\Phi_0$ and $\Phi$
we refer mainly to Triebel \cite[3.3]{Tr92}.\\\newline
The most prominent nontrivial examples (besides the one given in Remark
\ref{remunity}) of functions $\Phi_0$ and $\Phi$ satisfying
\eqref{condphi1} and \eqref{condphi2} are the classical local means. The name
comes from the compact support of $\Phi_0,\Phi$, which is admitted in the following statement.

\begin{corollary}\label{cor1} Let $p,q,s$ as in Theorem \ref{main1}. Let further
$k_0, k^0 \in \mathcal{S}(\R)$ such that
\be\label{f10}
    \cf k_0 (0), \cf k^0(0)\neq 0\,
\ee and define
\be\nonumber
  \Phi_0 = k_0\quad\mbox{and}\quad \Phi = \Delta^N k^0
\ee
with $N\in \N$ such that $2N>s$\,. Then \eqref{CL}, \eqref{CP},
\eqref{CT}, \eqref{DP} and \eqref{DL} characterize $F^s_{p,q}(\R)$\,.
\end{corollary}
\begin{corollary}\label{cor2} Let $p,q,s$ as in Theorem \ref{main1}. 
Let further $\varphi_0\in \mathcal{S}(\R)$ be a non-increasing radial function satisfying
$$
 \varphi_0 (0) \neq 0\quad\mbox{and}\quad D^{\bar{\alpha}}\varphi_0(0) = 0
$$
for $1\leq |\bar{\alpha}|_1 \leq R$, where $R+1 > s$\,.
Define $\varphi := \varphi_0(\cdot)-\varphi_0(2\cdot)$ and put $\Phi_0 := \cf^{-1}\varphi_0$ and 
$\Phi := \cf^{-1}\varphi$.
Then \eqref{DP} and \eqref{DL} characterize $F^s_{p,q}(\R)$\,.
\end{corollary}

\subsection{Proofs}
We give the proof for Theorem \ref{main1} in full detail. The proof
of Theorem \ref{main1b} is similar and even less technical. Let us 
refer to the respective paragraph for the necessary
modifications. The proofs in the Besov scale are analogous, so we
omit them completely. The proof technique is a modification of the one in Rychkov
\cite{Ry99a}, where he proved the discrete case, i.e., that
\eqref{DP} and \eqref{DL} characterize $F^s_{p,q}(\R)$. However,
Hansen \cite[Rem. 3.2.4]{HaDiss10} recently observed that the
arguments used for proving (34) in \cite{Ry99a} are
somehow problematic. The finiteness of the Peetre maximal function
is assumed which is not true in general under the stated
assumptions. Consider for instance in dimension $d=1$ the functions
$$
    \Psi_0(t) = \Psi_1(t) = e^{-t^2}\,
$$
and, if $a>0$ is given, the tempered distribution $f(t) = |t|^n$
with $a<n\in \mathbb{N}$. Then $(\Psi^{\ast}_kf)_a(x)$ is infinite in every point $x\in \re$.
The mentioned
incorrect argument was inherited to some
subsequent papers dealing with similar topics, for instance
\cite{Ba03}, \cite{Ke09} and \cite{Vy06}. Anyhow, the stated results
hold true. There is an alternative method to prove the crucial inequality (34) which avoids Lemma 3 in \cite{Ry99a}.
It is given in Rychkov \cite{Ry99b} as well as
\cite{Ry01}\,. A variant of this method, which is originally due to 
Str{\"o}mberg, Torchinsky \cite[Chapt. V]{StTo89}, is also used in
our proof below.\\

We start with a convolution type inequality which will be often needed
below. The following lemma is essentially Lemma 2 in \cite{Ry99a}.
\begin{lemma}\label{help2} Let $0<p,q\leq \infty$ and $\delta>0$. Let $\{g_{k}\}_{k\in
\n}$ be a sequence of non-negative measurable functions on $\R$ and
put
$$
    G_{\ell}(x) = \sum\limits_{k\in \zz}
    2^{-|k-\ell|\delta}g_{k}(x)\quad,\quad x\in
    \R\,,\ell \in \zz.
$$
Then there is some constant $C = C(p,q,\delta)$, such that
$$
     \|\{G_{\ell}\}_{\ell}|\ell_q(L_p(\R))\| \leq
     C\|\{g_{k}\}_{k}|\ell_q(L_p(\R))\|
$$
and
$$
     \|\{G_{\ell}\}_{\ell}|L_p(\ell_q,\R)\| \leq C\|\{g_{k}\}_{k}|L_p(\ell_q,\R)\|\
$$
hold true.
\end{lemma}

\subsubsection*{Proof of Theorem \ref{main1}}
To begin with we prove the equivalence of the 
characterizations \eqref{CL}, \eqref{CP}, \eqref{DP} and \eqref{DL}
for the same system $(\Phi_0,\Phi)$. The next step is to change 
from the system $(\Phi_0,\Phi)$ to a second one $(\Psi_0,\Psi)$ satisfying 
\eqref{condphi1}, \eqref{condphi2} 
within the characterization \eqref{DP}. 
The equivalence of \eqref{DP} and \eqref{DL} was the original proof
by Rychkov in \cite{Ry99a}. Since Definition \ref{inhom} can be seen as a special case
of \eqref{DL}, we have that
\eqref{CL}, \eqref{CP}, \eqref{DP} and \eqref{DL} generate the same
space for all pairs $(\Phi_0,\Phi)$ satisfying
\eqref{condphi1} and \eqref{condphi2}, namely $F^{s}_{p,q}(\R)$. 
It remains to prove that \eqref{CT} is equivalent to the rest. \frei\\\frei\\
{\em Step 1.}\\
We are going to prove the relations 
\begin{equation}\label{eq50}
   \|f|F^s_{p,q} \|_1 \asymp \|f|F^s_{p,q}\|_2 \asymp \|f|F^s_{p,q}\|_4 \asymp \|f|F^s_{p,q}\|_5
\end{equation}
for every $f\in \mathcal{S}'(\R)$\,. We just give the proof of $\|f|F^s_{p,q} \|_1 \asymp \|f|F^s_{p,q}\|_2$ in detail since the 
remaining equivalences are analogous. \\\newline
{\em Substep 1.1.}\\
 Put $\varphi_0 = \cf \Phi_0$ and $\varphi_{\ell} = (\cf
\Phi)(2^{-\ell}\cdot)$ if $\ell \geq 1$. Because of
\eqref{condphi1} it is possible to find functions $\psi_0,\psi \in
\mathcal{S}(\R)$ with $\supp \psi_0 \subset \{\xi\in \R~:~|\xi| \leq
2\varepsilon\}$, $\supp \psi \subset \{\xi\in \R~:~ \varepsilon/2
\leq |\xi| \leq 2\varepsilon\}$ and $\psi_\ell (x) = \psi(2^
{-\ell}x)$ such that
$$
    \sum\limits_{\ell \in \n} \varphi_{\ell}(\xi)\cdot \psi_{\ell}(\xi) = 1\,.
$$
We need a bit more. Fix a $1 \leq t \leq 2$. Clearly, we also have
$$
    \sum\limits_{\ell \in \n} \varphi_{\ell}(t\xi)\cdot \psi_{\ell}(t\xi) = 1\,
$$
for all $\xi \in \R$\,. With $\Psi_0 = \cf^{-1} \psi_0$ and $\Psi =
\cf^{-1}\psi$ we obtain then
\begin{equation}\nonumber
    g = \sum\limits_{m\in \n} (\Psi_{m})_t\ast
    (\Phi_{m})_t\ast g\,.
\end{equation}
We dilate this identity with $2^{\ell}$, i.e., $g_{\ell}(\eta) =
g(2^{-\ell d}\eta(2^{-\ell}\cdot))$ for $\eta \in \mathcal{S}(\R)$. An elementary calculation gives
\begin{equation}\label{f1}
    g_{\ell} = \sum\limits_{m\in \n} (\Psi_{m})_{t2^{-\ell}}\ast
    (\Phi_{m})_{t2^{-\ell}}\ast g_{\ell}\,
\end{equation}
for every $g \in \mathcal{S}'(\R)$. Obviously, we can rewrite
\eqref{f1} to obtain
\begin{equation}\label{f12}
    g  = \sum\limits_{m\in \n} (\Psi_{m})_{t2^{-\ell}}\ast
    (\Phi_{m})_{t2^{-\ell}}\ast g\,
\end{equation}
for all $g\in \mathcal{S}'(\R)$\,. Let us now choose $g = (\Phi_{\ell})_t \ast
f$ which gives for all $f\in \mathcal{S}'(\R)$ the identity
\begin{equation}\label{eq15}
    (\Phi_{\ell})_t \ast f  = \sum\limits_{m\in \n} (\Phi_\ell)_t \ast (\Psi_{m})_{t2^{-\ell}}\ast
    (\Phi_{m})_{t2^{-\ell}}\ast f\,.
\end{equation}
For $m,\ell\in \n$ we define
\begin{equation}\label{f13}
    \Lambda_{m,\ell}(x) = \left\{\begin{array}{lcl}
    2^{\ell d}\Phi_0(2^{\ell}x)&:&m=0\\
    \Phi_{\ell}(x)&:& m>0
    \end{array}\right.\quad,\quad x\in \R\,.
\end{equation}
Clearly, we have
\begin{equation}\nonumber
  (\Phi_{\ell})_t\ast (\Phi_m)_{t2^{-\ell}}
  =(\Lambda_{m,\ell})_t \ast (\Phi_{m+\ell})_t.
\end{equation}
Plugging this into \eqref{eq15} we end up with the pointwise
representation ($\ell \in \N$)
\begin{equation}
  \label{eq28}
  \begin{split}
    ((\Phi_{\ell})_t\ast f)(y) &= \sum\limits_{m\in \n} ((\Psi_m)_{2^{-\ell}t}
    \ast (\Lambda_{m,\ell})_t \ast
    (\Phi_{m+\ell})_t\ast f)(y)\\
    &=\sum\limits_{m\in \n} [(\Psi_m)_{2^{-\ell}t}
    \ast (\Lambda_{m,\ell})_t] \ast
    ((\Phi_{m+\ell})_t\ast f)(y)\\
    &=\sum\limits_{m\in \n}\int\limits_{\R} [(\Psi_m)_{2^{-\ell}t}\ast
    (\Lambda_{m,\ell})_t](y-z)\cdot ((\Phi_{m+\ell})_t\ast
    f)(z)\,dz
  \end{split}
\end{equation}
for all $y\in \R$\,. Let us mention that the case $\ell = 0$ plays a particular role. In this case we have to replace $(\Phi_{\ell})_t$ 
by $\Phi_{\ell}$ in \eqref{eq15} and \eqref{eq28}, $(\Phi_{m+\ell})_t$ by $\Phi_{m+\ell}$ in \eqref{eq28} if $m=0$ and finally $(\Lambda_{m,\ell})_t$ by
$\Lambda_{m,\ell}$ if $m>0$.  \frei\\\frei\\
{\em Substep 1.2.} Let us prove the following important inequality
first. For every $r>0$ and every $N \in \n$ we have
\begin{equation}\label{eq16}
    |((\Phi_{\ell})_t\ast f)(x)|^r \leq
    c\sum\limits_{k\in \n} 2^{-kNr}2^{(k+\ell)d}
    \int\limits_{\R}\frac{
    |((\Phi_{k+\ell})_t\ast
    f)(y)|^r}{(1+2^{\ell}|x-y|)^{Nr}}\,dy\,,
\end{equation}
where $c$ is independent of $f \in \mathcal{S}'(\R)$, $x\in \R$ and $\ell\in
\n$\,. Again the case $\ell = 0$ has to be treated separately according to the remark after \eqref{eq28}. The representation \eqref{eq28} will be the starting
point to prove \eqref{eq16}. Namely, we have for $y\in \R$
\begin{eqnarray}
        \nonumber
        |((\Phi_{\ell})_t\ast f)(y)| &\leq&
        \sum\limits_{m\in \n}\int\limits_{\R} |((\Psi_{m})_{t2^{-\ell}}\ast
        (\Lambda_{m,\ell})_t)(y-z)|\cdot |((\Phi_{m+\ell})_t\ast
        f)(z)|\,dz\nonumber\\
        &\leq& \sum\limits_{m\in \n} S_{m,\ell,t}\int\limits_{\R} \frac{|((\Phi_{m+\ell})_t\ast
        f)(z)|}{(1+2^{\ell}|y-z|)^N}\,dz\,,
        \label{eq18}
\end{eqnarray}
where
$$
  S_{m,\ell,t} = \sup\limits_{x\in \R} |[(\Psi_{m})_{2^{-\ell}t}\ast
        (\Lambda_{m,\ell})_t](x)|\cdot
        (1+2^{\ell}|x|)^{N}\,.
$$
Elementary properties of the convolution yield (compare with
\eqref{eq8})
\begin{equation}
\nonumber
 \begin{split}
    S_{m,\ell,t} &= \frac{2^{\ell d}}{t^d}\sup\limits_{x\in \R} |(\Psi_{m} \ast (\Lambda_{m,\ell})_{2^{\ell}})(x2^{\ell}/t)|\cdot
    (1+2^{\ell}|x|)^{N}\\
    &= \frac{2^{\ell d}}{t^d}\sup\limits_{x\in \R} |(\Psi_{m} \ast \eta_{m,\ell})(x)|\cdot
    (1+|tx|)^{N}\,,
 \end{split}
\end{equation}
where
$$
    \eta_{m,\ell}(x) = \left\{\begin{array}{lcl}
        \Phi(x)&:&m>0, \ell>0\,,\\
        \Phi_0(x)&:&\mbox{otherwise}\,.
    \end{array}\right.
$$
With Lemma \ref{help1} we see
$$
    S_{m,\ell,t} \leq  C_{N}2^{\ell d}2^{-mN}
$$
and put it into \eqref{eq18} to obtain
\begin{equation}
   \label{eq29}
   |((\Phi_{\ell})_t\ast f)(y)|\leq C_{N}
   \sum\limits_{m\in \n}2^{-mN}\int\limits_{\R}
   \frac{2^{(m+\ell)d}|((\Phi_{m+\ell})_t\ast
   f)(z)|}{(1+2^{\ell}|y-z|)^{N}}\,dz\,
\end{equation}
with the appropriate modification in case $\ell = 0$. To continue we prefer the strategy used by Rychkov in \cite[Thm. 3.2]{Ry99b} and
\cite[Lem. 2.9]{Ry01}.\\
Let us replace $\ell$ by $k+\ell$ in \eqref{eq29} and
multiply on both sides with $2^{-kN}$. Then we can estimate
\begin{eqnarray}
    \label{eq27}
    \hspace{-1cm}2^{-kN}|((\Phi_{k+\ell})_t\ast f)(y)|
    &\leq& C_{N}\sum\limits_{m\in \n}
    2^{-kN}2^{-mN}
    \int\limits_{\R} \frac{2^{(m+k+\ell)d}|((\Phi_{m+k+\ell})_t\ast
   f)(z)|}{(1+2^{k+\ell}|y-z|)^{N}}\,dz\\
   &\leq& C_{N}\sum\limits_{m\in \n}2^{-(m+k)N}
   \int\limits_{\R} \frac{2^{(m+k+\ell)d}|((\Phi_{m+k+\ell})_t\ast
   f)(z)|}{(1+2^{\ell}|y-z|)^{N}}\,dz\nonumber\\
   &=& C_{N}\sum\limits_{m\in k+\n}2^{-mN}
   \int\limits_{\R} \frac{2^{(m+\ell)d}|((\Phi_{m+\ell})_t\ast
   f)(z)|}{(1+2^{\ell}|y-z|)^{N}}\,dz\nonumber\\
   &\leq&
   C_{N}\sum\limits_{m\in \n}2^{-mN}
   \int\limits_{\R} \frac{2^{(m+\ell)d}|((\Phi_{m+\ell})_t\ast
   f)(z)|}{(1+2^{\ell}|y-z|)^{N}}\,dz\,.
   \label{eq21}
\end{eqnarray}
Next, we apply the elementary inequalities
\begin{equation}\label{eq19}
    (1+2^{\ell}|y-z|)\cdot (1+2^{\ell}|x-y|)\geq
    (1+2^{\ell}|x-z|)\,,
\end{equation}
\begin{equation}
    \nonumber
    \begin{split}
      |((\Phi_{m+\ell})_t \ast f)(z)| \leq&
      |((\Phi_{m+\ell})_t\ast f)(z)|^r (1+2^{\ell}|x-z|)^{N(1-r)}\\
      &~~\times\sup\limits_{y\in \R}\frac{|((\Phi_{m+\ell})_t\ast f)(y)|^{1-r}}{
      (1+2^{\ell}|x-y|)^{N(1-r)}}\,,
    \end{split}
\end{equation}
where $0<r\leq 1$. Let us define the maximal function
\begin{equation}\label{Str}
    M_{\ell,N}(x,t) = \sup\limits_{k\in
    \n}\sup\limits_{y\in \R} 2^{-kN}
    \frac{|((\Phi_{k+\ell})_t\ast
    f)(y)|}{(1+2^{\ell}|x-y|)^{N}}\quad,\quad
    x\in \R\,,
\end{equation}
and estimate
\begin{eqnarray}
     \label{eq26}
     M_{\ell,N}(x,t) &\leq& C_{N}\sum\limits_{m\in \n}2^{-mN}
     \int\limits_{\R} \frac{2^{(m+\ell)d}|((\Phi_{m+\ell})_t\ast
     f)(z)|}{(1+2^{\ell}|x-z|)^{N}}\,dz\\
     &\leq& C_{N}\sum\limits_{m\in
     \n}2^{-mNr}\Big(2^{-mN}\sup\limits_{y\in \R}
     \frac{|((\Phi_{m+\ell})_t\ast f)(y)|}{
     (1+2^{\ell}|x-y|)^{N}}\Big)^{1-r}\label{eq22}\\
     &~~~&\int\limits_{\R} \frac{2^{(m+\ell)d}|((\Phi_{m+\ell})_t\ast
     f)(z)|^r}{(1+2^{\ell}|x-z|)^{Nr}}\,dz\,.\nonumber
\end{eqnarray}
Observe that we can estimate the term $(...)^{1-r}$ in the
right-hand side of \eqref{eq22} by $M_{\ell,N}(x,t)^{1-r}$. Hence,
if $M_{\ell,N}(x,t)<\infty$ we obtain from \eqref{eq22}
\begin{equation}\label{eq25}
    M_{\ell,N}(x,t)^r \leq C_{N}\sum\limits_{m\in
    \n}2^{-mNr}\int\limits_{\R} \frac{2^{(m+\ell)d}|((\Phi_{m+\ell})_t\ast
    f)(z)|^r}{(1+2^{\ell}|x-z|)^{Nr}}\,dz\,,
\end{equation}
where $C_{N}$ is independent of $x$, $f$, $\ell$ and $t \in [1,2]$. We
claim that there exists $N^f \in \n$ such that $M_{\ell,N}(x,t) <
\infty$ for all $N\geq N^f$\,. Indeed, we use that $f\in \mathcal{S}'(\R)$,
i.e., there is an $M\in \n$ and $c_f>0$ such that
$$
    |((\Phi_{k+\ell})_t\ast f)(y)| \leq c_f\sup\limits_{|\bar{\alpha}|_1\leq
    M}\sup\limits_{z\in \R}
    |D^{\bar{\alpha}}\Phi_{k+\ell}(z)|\cdot(1+|y-z|)^M\,,
$$
see \eqref{eq100}.
Assuming $N>M$ we estimate as follows
\begin{eqnarray}
      |((\Phi_{\ell})_t\ast f)(x)| &\leq&
      M_{\ell,N}(x,t)\label{eq24}\\
      &\leq&
      c\sup\limits_{k\in \n}\sup\limits_{y\in \R} 2^{-kN}
      \frac{|((\Phi_{k+\ell})_t\ast
      f)(y)|}{(1+|x-y|^2)^{\frac{N}{2}}}\nonumber\\
      &\leq& c\sup\limits_{k\in \n}\sup\limits_{y\in \R}
      2^{-kN}2^{(k+\ell)(M+d)}\sup\limits_{z\in \R}\sup\limits_{|\bar{\alpha}|_1\leq M}
      \frac{|D^{\bar{\alpha}}\gamma_{k+\ell}(z)|\cdot
      (1+|y-z|)^M}{(1+|x-y|)^{N
      }}\nonumber\\
      &\leq& c2^{\ell(M+d)}\sup\limits_{k\in \n}
      \sup\limits_{z\in \R}\sup\limits_{|\bar{\alpha}|_1\leq M}
      |D^{\bar{\alpha}}\gamma_{k+\ell}(z)|(1+|x-z|)^{N}\,,\nonumber
\end{eqnarray}
where we again used the inequality (compare with \eqref{eq19})
$$
    1+|y-z| \leq (1+|x-y|)(1+|x-z|)
$$
and put
$$
    \gamma_{\ell}(t) = \left\{\begin{array}{lcl}
        \Phi_0(t)&:&\ell=0\\
        \Phi(t)&:&\ell>0
    \end{array}\right..
$$
Hence $\gamma_{k+\ell}$ gives us only two different functions from
$\mathcal{S}(\R)$. This implies the boundedness of $M_{\ell,N}(x,t)$ for $x\in
\R$ if $N>M = N^f$. Therefore, \eqref{eq25} together with
\eqref{eq24} yield \eqref{eq16} with $c=C_{N}$, independent of $x$,
$f$ and $\ell$, for all $N \geq N^f$. But this is not yet what we
want. Observe that the right-hand side of \eqref{eq16} decreases as
$N$ increases. Therefore, we have \eqref{eq16} \underline{for all}
$N \in \n$ but with $c = c(f) = C_{N^f}$ \underline{depending on
$f$}. This is still not yet what we want. Now we argue as follows:
Starting with \eqref{eq16} where $c=c(f)$ and $N\in \n$ arbitrary,
we apply the same arguments as used from \eqref{eq27} to
\eqref{eq21}, switch to the maximal function \eqref{Str} with the
help of \eqref{eq19} and finish with \eqref{eq25} instead of
\eqref{eq26} but with a constant that depends on $f$. But this does
not matter now. Important is, that a finite right-hand side of
\eqref{eq25} (which is the same as
rhs\eqref{eq16}) implies $M_{\ell,N}(x,t)<\infty$.\\
We assume $\mbox{rhs}\eqref{eq16}<\infty$. Otherwise there is
nothing to prove in \eqref{eq16}. Returning to \eqref{eq22} and
having in mind that now $M_{\ell,N}(x,t)<\infty$, we end up with
\eqref{eq25} for all $N$ and $C_{N}$ independent of $f$. Finally,
from \eqref{eq25} we obtain \eqref{eq16} and are done in case
$0<r\leq 1$.\frei\\ Of course, \eqref{eq16} also holds true for
$r>1$ with a much simpler proof. In that case, we use \eqref{eq29}
with $N+d+\varepsilon$ instead of $N$ and apply H{\"o}lder's inequality with
respect to $1/r+1/r' = 1$ first for integrals and then for
sums.\\\newline {\em Substep 1.3.}\\ The inequality \eqref{eq16}
implies immediately a stronger version of itself. Using \eqref{eq19}
again we obtain for $a\leq N$ and $\ell \in \N$
\begin{equation}\label{f2}
    (\Phi^{\ast}_{2^{-\ell}t}f)_{a}(x)^r\leq
    c\sum\limits_{k\in \n} 2^{-kNr}2^{(k+\ell)d}
    \int\limits_{\R}\frac{
    |((\Phi_{k+\ell})_t\ast
    f)(y)|^r}{(1+2^{\ell}|x-y|)^{a
    r}}\,dy\,.
\end{equation}
In case $\ell = 0$ we have to replace $(\Phi^{\ast}_{2^{-\ell}t}f)_{a}(x)$ by $(\Phi^{\ast}_0f)_a(x)$ on the left-hand side and
$(\Phi_{k+\ell})_t$ by $\Phi_{k+\ell}=\Phi_0$ for $k=0$ on the right-hand side. We proved, that the inequality \eqref{f2} holds for all $t \in
[1,2]$ where $c>0$ is independent of $t$. If we choose
$r<\min\{p,q\}$, we can apply the norm
$$
    \Big(\int\limits_{1}^2 |\cdot|^{q/r}\frac{dt}{t}\Big)^{r/q}\,.
$$
on both sides and use Minkowski's inequality for integrals, which
yields
\begin{equation}\label{f3}
    \Big(\int\limits_{1}^2|(\Phi^{\ast}_{2^{-\ell}t}f)_{a}(x)|^q\,\frac{dt}{t}\Big)^{r/q} \leq
    c\sum\limits_{k\in \n} 2^{-kNr}2^{(k+\ell)d}
    \int\limits_{\R}\frac{
    \left(\int\limits_1^2|((\Phi_{k+\ell})_t\ast
    f)(y)|^q\,\frac{dt}{t}\right)^{r/q}}{(1+2^{\ell}|x-y|)^{a
    r}}\,dy\,.
\end{equation}
If $a r>d$ then we have
$$
        g_{\ell}(y) = \frac{2^{d\ell}}{{(1+2^{\ell}|y|)^{ar}}}
        \in L_1(\R)\,
$$
and we observe
\begin{equation}
  \nonumber
  \begin{split}
    \Big(\int\limits_{1}^2|2^{\ell s}(\Phi^{\ast}_{2^{-\ell}t}f)_{a}(x)|^q\,\frac{dt}{t}\Big)^{r/q}
    \leq
    c\sum\limits_{k\in \n} 2^{-kNr}2^{kd}2^{\ell s r}
    \Big[g_{\ell} \ast
    \Big(\int\limits_{1}^2|2^{\ell s}((\Phi_{k+\ell})_t\ast f)(\cdot)|^q\,\frac{dt}{t}\Big)^{r/q}\Big](x)\,.
 \end{split}
 \end{equation}   
Now we use a well-known majorant property in order to estimate the convolution on the right-hand side 
by the Hardy-Littlewood maximal function 
(see Paragraph \ref{sectmax} and
\cite[Chapt. 2]{StWe71}). This yields
\vspace{-0.3cm}
\begin{equation}
  \nonumber  
    \Big(\int\limits_{1}^2|2^{\ell s}(\Phi^{\ast}_{2^{-\ell}t}f)_{a}(x)|^q\,\frac{dt}{t}\Big)^{r/q}
    \leq c\sum\limits_{k\in \n} 2^{\ell r s}2^{k(-Nr+d)}
    M\Big[\Big(\int\limits_{1}^2|((\Phi_{k+\ell})_t\ast f)(\cdot)|^q\,\frac{dt}{t}\Big)^{r/q}\Big](x)\,.
 \end{equation}   
An index shift on the right-hand side gives
\begin{equation}
  \nonumber
  \begin{split}
    &\Big(\int\limits_{1}^2|2^{\ell s}(\Phi^{\ast}_{2^{-\ell}t}f)_{a}(x)|^q\,\frac{dt}{t}\Big)^{r/q}\\
    &\hspace{1cm}\leq c\sum\limits_{k\in \ell+\n} 2^{\ell r s}2^{(k-\ell)(-Nr+d)}
    M\Big[\Big(\int\limits_{1}^2|((\Phi_{k})_t \ast f)(\cdot)|^q\,\frac{dt}{t}\Big)^{r/q}\Big](x)\\
    &\hspace{1cm}=c\sum\limits_{k\in \ell+\n} 2^{(\ell-k)(Nr-d+rs)}2^{krs}
    M\Big[\Big(\int\limits_{1}^2|((\Phi_{k})_t \ast f)(\cdot)|^q\,\frac{dt}{t}\Big)^{r/q}\Big](x)\,.
 \end{split}
\end{equation}
Choose now $d/a < r < \min\{p,q\}$, $N>\max\{0,-s\}+a$ and put
$$
    \delta = N+s-d/r>0\,.
$$
We obtain for $\ell \in \N$
$$
  \Big(\int\limits_{1}^2|2^{\ell s}(\Phi^{\ast}_{2^{-\ell}t}f)_{a}(x)|^q\,\frac{dt}{t}\Big)^{r/q}
   \leq c\sum\limits_{k\in \N} 2^{-\delta r|\ell -k|}2^{krs}M\Big[\Big(\int\limits_{1}^2|((\Phi_{k})_t \ast f)(\cdot)|^q\,\frac{dt}{t}\Big)^{r/q}\Big](x)\,.
$$
Now we apply Lemma \ref{help2} in $L_{p/r}(\ell_{q/r},\R)$ which
yields
\begin{equation}\nonumber
   \begin{split}
    &\Big\|\Big(\int\limits_{1}^2|2^{\ell s}(\Phi^{\ast}_{2^{-\ell}t}f)_{a}(x)|^q\,\frac{dt}{t}\Big)^{r/q}|L_{p/r}(\ell_{q/r})\Big\|\\
    &\hspace{1cm}\leq c \Big\|M\Big[\Big(\int\limits_{1}^2
    |2^{ks}((\Phi_{k})_t \ast f)(\cdot)|^q\,\frac{dt}{t}\Big)^{r/q}\Big]|L_{p/r}(\ell_{q/r})\Big\|\,.
   \end{split}
\end{equation}
The Fefferman-Stein inequality (see Paragraph
\ref{sectmax}/Theorem \ref{feffstein}, having in mind that
$p/r,q/r>1$) gives
\begin{equation}
\nonumber
   \begin{split}
    &\Big\|\Big(\int\limits_{1}^2|2^{\ell s}(\Phi^{\ast}_{2^{-\ell}t}f)_{a}(x)|^q\,\frac{dt}{t}\Big)^{1/q}|L_{p}(\ell_{q})\Big\|^r\\
    &\hspace{1cm}\lesssim\Big\|M\Big[\Big(\int\limits_{1}^2|2^{\ell s}(\Phi^{\ast}_{2^{-\ell}t}f)_{a}(x)|^q\,\frac{dt}{t}\Big)^{r/q}\Big]|L_{p/r}(\ell_{q/r})\Big\|\\
    &\hspace{1cm}\lesssim \Big\|\Big(\int\limits_{1}^2
    |2^{ks}((\Phi_{k})_t \ast f)(\cdot)|^q\,\frac{dt}{t}\Big)^{r/q}|L_{p/r}(\ell_{q/r})\Big\|\\
    &\hspace{1cm} = \Big\|\Big(\int\limits_{1}^2
    |2^{ks}((\Phi_{k})_t \ast f)(\cdot)|^q\,\frac{dt}{t}\Big)^{1/q}|L_{p}(\ell_{q})\Big\|^r\,.
   \end{split}
\end{equation}
Hence, we obtain
\begin{equation}\nonumber
  \begin{split}
      \Big\|\Big(\int\limits_{0}^1 |\lambda^{-sq}(\Phi^{\ast}_{\lambda} f)_{a}(x)|^q\,\frac{d \lambda}{\lambda}\Big)^{1/q}|L_p(\R)\Big\|
      &\asymp\Big\|\Big(\sum\limits_{\ell=1}^{\infty}\int\limits_{1}^2|2^{\ell s}(\Phi^{\ast}_{2^{-\ell}t}f)_{a}(x)|^q\,\frac{dt}{t}
      \Big)^{1/q}|L_{p})\Big\|\\
      &=\Big\|\Big(\int\limits_{1}^2|2^{\ell s}(\Phi^{\ast}_{2^{-\ell}t}f)_{a}(x)|^q\,\frac{dt}{t}\Big)^{1/q}|L_{p}(\ell_{q})\Big\|\\
      &\lesssim \Big\|\Big(\int\limits_{0}^1 |\lambda^{-sq}(\Phi_{\lambda} \ast f)(x)|^q\,\frac{d \lambda}{\lambda}\Big)^{1/q}|L_p(\R)\Big\|\,.
  \end{split}
\end{equation}
The summand $\|(\Phi^{\ast}_0f)_a|L_p(\R)\|$ can be estimated similar using \eqref{f2} in case $\ell = 0$. This proves $\|f|F^s_{p,q}(\R)\|_2 \lesssim \|f|F^s_{p,q}(\R)\|_1$. With slight modifications of the argument we prove as well
$\|f|F^s_{p,q}(\R)\|_2 \lesssim \|f|F^s_{p,q}(\R)\|_5$, $\|f|F^s_{p,q}(\R)\|_4 \lesssim \|f|F^s_{p,q}(\R)\|_1$, and $\|f|F^s_{p,q}(\R)\|_4 \lesssim \|f|F^s_{p,q}(\R)\|_5$. The inequalities $\|f|F^s_{p,q}(\R)\|_5 \lesssim \|f|F^s_{p,q}(\R)\|_4$ and $\|f|F^s_{p,q}(\R)\|_1 \lesssim \|f|F^s_{p,q}(\R)\|_2$  
are immediate. This finishes the proof of \eqref{eq50}.\\\frei\\
{\em Step 2.} Let $\Psi_0, \Psi \in \mathcal{S}(\R)$ be functions satisfying \eqref{condphi2}. 
Indeed, we do not need \eqref{condphi1} for the following inequality 
\begin{equation}\label{eq51}
    \|f|F^s_{p,q}(\R)\|^{\Psi}_4 \lesssim \|f|F^s_{p,q}(\R)\|^{\Phi}_4  \,
\end{equation}
which holds true for all $f\in \mathcal{S}'(\R)$. We decompose $f$ similar as in Step 1.
Exploiting the property \eqref{condphi1} for the system $(\Phi_0,\Phi)$ we
find $\mathcal{S}(\R)$-functions $\lambda_0, \lambda \in \mathcal{S}(\R)$ such that
$\supp \lambda_0 \subset \{\xi \in \R~:~|\xi| \leq 2\varepsilon\}$
and $\supp \lambda \subset \{\xi \in \R~:~\varepsilon /2 \leq |\xi|
\leq 2\varepsilon\}$ and
$$
    \sum\limits_{k\in \n} \lambda_{k}(\xi)\cdot
    \varphi_{k}(\xi) = 1\,
$$
for $\xi\in \R$\,. Putting $\Lambda_0 = \cf^{-1}
\lambda_0$ and $\Lambda = \cf^{-1} \lambda$ we obtain the decomposition
\be\label{eq14}
    g = \sum\limits_{k\in \n} \Lambda_{k} \ast
    \Phi_{k}\ast g
\ee 
for every $g \in \mathcal{S}'(\R)$\,. We put $g = \Psi_\ell \ast f$ for $\ell \in \n$ and
see
\be\label{f4}
    \Psi_{\ell} \ast f = \sum\limits_{k \in \n} \Psi_{\ell}\ast
    \Lambda_{k}\ast \Phi_{k}\ast f\,.
\ee
Now we estimate as follows 
\begin{equation}
\nonumber
    \begin{split}
       |((\Psi_{\ell}\ast \Lambda_{k})\ast (\Phi_k \ast f))(y)|
       &\leq
       \int\limits_{\R}|(\Psi_{\ell}\ast\Lambda_{k})(z)\cdot |(\Phi_k \ast f)(y-z)|\,dz\\
       &\leq (\Phi^{\ast}_{k}f)_a(y)
       \int\limits_{\R}|(\Phi_{\ell}\ast\Lambda_{k})(z)|\cdot(1+2^{k}|z|)^{a}\,dz\\
       &\leq (\Phi^{\ast}_{k}f)_a(y) J_{\ell,k\,,}
    \end{split}
\end{equation}
where
\begin{equation}\nonumber
   J_{\ell,k} = \int\limits_{\R} |(\Psi_{\ell}\ast
   \Lambda_{k})(z)|(1+2^{k}|z|)^{a}\,dz\,.
\end{equation}
We first observe that for $x\in \R$ and functions $\mu,\eta \in
\mathcal{S}(\R)$ the following identity holds true for $u,v>0$

\be\label{eq8}
  (\mu_u \ast \eta_v)(z) = \frac{1}{u^d} [\mu \ast \eta_{v/u}](z/u)
  = \frac{1}{v^d}[\mu_{u/v}\ast \eta](z/v)\,.
\ee This yields in case $\ell\geq k$ (with a minor change if $k = 0$) \be\nonumber
  \begin{split}
    J_{\ell,k} &= \int\limits_{\R}|[\Psi_{\ell-k} \ast
    \Lambda](z)|(1+|z|)^a\\
    &\lesssim \sup\limits_{z\in \R} |[\Psi_{\ell-k}\ast
    \Lambda](z)|(1+|z|)^{a+d+1}\\
    &\lesssim 2^{(k-\ell)(R+1)}\,,
  \end{split}
\ee where we used Lemma \ref{help1} for the last estimate.\\
If $k > \ell$ we change the roles of $\Psi$ and $\Lambda$ to obtain
again with Lemma \ref{help1} (minor change if $\ell = 0$)
\begin{equation}
  \nonumber
  \begin{split}
    J_{\ell,k} &= \int\limits_{\R} |[\Psi\ast \Lambda_{k-\ell}](z)|
    (1+|2^{k-\ell}z|)^{a}\,dx\\
    &\lesssim 2^{(k-\ell)a}\sup\limits_{z\in \R} |[\Psi\ast \Lambda_{k-\ell}](z)|(1+|z|)^{a+d+1}\\
    &\lesssim 2^{(\ell-k)(L+1-a)}\,,
  \end{split}
\end{equation}
where $L$ can be chosen arbitrary large since $\Lambda$ satisfies
$(M_L)$ for every $L \in \N$ according to its construction. Let us
further use the estimate \be\nonumber
  \begin{split}
    (\Phi^{\ast}_{k}f)_a(y) &\leq
    (\Phi^{\ast}_{k}f)_a(x)(1+2^{k}|x-y|)^a\\
    &\lesssim
    (\Phi^{\ast}_{k}f)_a(x)(1+2^{\ell}|x-y|)^a\max\{1,2^{(k-\ell)a}\}\,.
  \end{split}
\ee Consequently,
\begin{equation}
    \nonumber
    \begin{split}
        \sup\limits_{y \in \R}\frac{2^{s\ell}|(\Psi_{\ell}\ast \Lambda_{k}\ast
        (\Phi_k\ast f))(y)|}{(1+2^{\ell}|x-y|)^a} &\lesssim
        2^{ks}(\Phi^{\ast}_{k}f)_a(x) 2^{(\ell-k) s}\max\{1,2^{(k-\ell)a}\}J_{\ell,k}\\
        &\leq 2^{ks}(\Phi^{\ast}_{k}f)_a(x)
        \left\{\begin{array}{lcl}
            2^{(\ell-k)(L+1-a+s)}&:&k>\ell\\
            2^{(k-\ell)(R+1-s)}&:&\ell \geq k
        \end{array}\right..
    \end{split}
\end{equation}
Plugging this into \eqref{f4}, choosing $L \geq a+|s|$ and $\delta =
\min\{1,R+1-s\}$ we obtain the inequality 
\begin{equation}\label{eq9}
   2^{\ell s}(\Psi_{\ell}^{\ast}f)_a(x) \lesssim \sum\limits_{k = 0}^{\infty}2^{-|k-\ell|\delta}
   2^{ks}(\Phi^{\ast}_kf)_a(x)\,.
\end{equation}
for all $x\in \R$\,. Applying Lemma \ref{help2} gives \eqref{eq51}\,.\frei\\\frei\\
{\em Step 3. } What remains is to show that \eqref{CT} is
equivalent to the rest.\frei\\\frei\\
{\em Substep 3.1. } Let us prove
\be\label{f5}
    \|f|F^s_{p,q}(\R)\|_2 \lesssim \|f|F^s_{p,q}(\R)\|_3\,.
\ee
We return to \eqref{f2} in Substep 1.3. If
$|z|<2^{-(\ell+k)}t$ formula \eqref{f2} implies by shift in the
integral the following
\begin{equation}\label{eq-33}
    (\Phi^{\ast}_{2^{-\ell}t}f)_{a}(x)^r\leq
    C_N\sum\limits_{k\in \n} 2^{-k(N-a)r}2^{(k+\ell)d}
    \int\limits_{\R}\frac{
    |((\Phi_{k+\ell})_t\ast
    f)(y+z)|^r}{(1+2^{\ell}|x-y|)^{a
    r}}\,dy\,.
\end{equation}
Indeed, we have
\begin{equation}\nonumber
  \begin{split}
    1+2^{\ell}|x-y| &\leq 1+2^{\ell}(|x-(y+z)|+|z|)\\
    &\lesssim 1+2^{\ell}(|x-(y+z)|)+2^{-k}\\
    &\lesssim 1+2^{\ell}(|x-(y+z)|)\,.
  \end{split}
\end{equation}
Where the last estimate follows from the fact that $k\in \N_0$ in
the sum. Instead of the integral $(\int_{1}^2
|\cdot|^{q/r}\,dt/t)^{r/q}$ we now take on both sides of \eqref{eq-33}
the norm
$$
    \Big(\int\limits_{1}^2 \int\limits_{|z|<t} |\cdot|^{q/r}\,dz
    \frac{dt}{t^{d+1}}\Big)^{r/q}\,.
$$
The integration over $z$ does not influence the left-hand side.
Instead of \eqref{f3} we obtain \be\nonumber
\begin{split}
    &\Big(\int\limits_{1}^2|(\Phi^{\ast}_{2^{-\ell}t}f)_{a}(x)|^q\,\frac{dt}{t}\Big)^{r/q}\\
    &\hspace{1cm}\leq c\sum\limits_{k\in \n} 2^{-kNs}2^{(k+\ell)d}
    \int\limits_{\R}\frac{
    \left(\int\limits_1^2 \int \limits_{|z|<t}|((\Phi_{k+\ell})_t\ast
    f)(y)|^q\,dz\,\frac{dt}{t^{d+1}}\right)^{r/q}}{(1+2^{\ell}|x-y|)^{a
   r}}\,dy\,.
\end{split}
\ee
We continue with analogous arguments as after \eqref{f3} and end up
with \eqref{f5}\,.\frei\\\frei\\
{\em Substep 3.2.} We prove $\|f|F^{s}_{p,q}(\R)\|_3 \lesssim
\|f|F^s_{p,q}(\R)\|_2$\,. Indeed, it is easy to see, that we have
for all $t>0$
\be\nonumber
  \begin{split}
    \frac{1}{t^d}\int\limits_{|z|<t} |(\Phi_t\ast f)(x+z)|\,dz
    &\lesssim \sup\limits_{|z|<t} \frac{|(\Phi_t\ast
    f)(x+z)|}{(1+|z|/t)^a}\\
    &\lesssim (\Phi_t^{\ast} f)_a(x)\,,
  \end{split}
\ee
and we are done. The proof is complete \eproof

\subsubsection*{Proof of Theorem \ref{main1b}}

The proof of Theorem \ref{main1b} is almost the same as the previous one. 
It is less technical since we do not have to
deal with a separate function $\Phi_0$ which causes several difficulties. 
However, there are still some technical obstacles
which have to be discussed.\\\newline 1. Although we are in the
homogeneous world, we use the same decomposition as used in
\eqref{f12}, even with the inhomogeneity $\Phi_0$. In the definition
of $\Lambda_{m,\ell}(x)$ in \eqref{f13} we have to put in addition $\Phi(x)$, if
$\ell=0$ and $m>0$. The consequence is equation \eqref{eq28} for
every $\ell\in \zz$. Hence, the inhomogeneity
is shifted to $\Lambda_{m,\ell}$\,. This yields \eqref{f2} for all
$\ell \in \zz$, where $k$ still runs through $\n$. We need this for
the argument in Substep 3.1. \\\newline 2. In contrast to the
previous decomposition, we use \eqref{eq14}, \eqref{f4} now for $k,\ell\in\zz$,
where $\Phi_0 = \Phi$ and $\Lambda_0 = \Lambda$. This
works since we assume $g\in \mathcal{S}_0'(\R)$\,. Now we can even prove 
$\|f|\dot{F}^s_{p,q}(\R)\|^{\Psi}_4 \lesssim \|f|\dot{F}^s_{p,q}(\R)\|^{\Phi}_2$
and vice versa. 
\eproof

\subsubsection*{Proof of Corollary \ref{cor1} and \ref{cor2}}

1. The proof of Corollary \ref{cor1} is immediate. We know that
$\Delta^N$ gives $(\sum_{k=1}^d |\xi_k|^2)^N$ as factor on the Fourier
side. This gives \eqref{condphi2} immediately and together with
\eqref{f10} we have \eqref{condphi1} for $\varepsilon>0$ small
enough.\\\newline 2. In the case of Corollary \ref{cor2} the
situation is a bit more involved. Clearly, Condition \eqref{condphi2}
holds true. But the problem here is, that \eqref{condphi1} may be
violated for all $\varepsilon>0$. However, we argue as follows. In
Step 2 in the proof above we have seen, that we do not
need \eqref{condphi1} for the system $(\Psi_0,\Psi)$. Hence, we can estimate
\eqref{DP} and \eqref{DL} from above by a further characterization of
$F^s_{p,q}(\R)$\,. For the remaining estimates we apply Theorem
\ref{main1} with the system $(\Phi_0, \tilde{\Phi})$ where
$$
    \tilde{\Phi} = \Phi_0(x) - \frac{1}{2^{kd}}\Phi_0(x/2^{k})\,,
$$
and $k\in \N$ is chosen in such a way that \eqref{condphi1} is
satisfied. What remains is a consequence of the fact that
$$
    \tilde{\Phi} = \Phi + \Phi_{-1} + ... + \Phi_{-(k-1)}\,.
$$
This type of argument is due to Triebel \cite[3.3.3]{Tr92}\,.\eproof

\section{Classical coorbit space theory}
\label{Classcoo} In \cite{FeGr86,FeGr89a,FeGr89b,Gr91} a general theory
of Banach spaces related to integrable group representations has been
developed. The ingredients are a locally compact group
$\mathcal{G}$ with identity $e$, a Hilbert space $\mathcal{H}$ and an irreducible,
unitary and continuous representation $\pi:\mathcal{G} \to \mathcal{L}(\mathcal{H})$, which
is at least integrable. One can associate a Banach space $\CoY$ to any solid,
translation-invariant Banach space $Y$ of functions on the group
$\mathcal{G}$. The main achievement
of this abstract theory is a powerful discretization machinery for $\CoY$, i.e., a universal approach to 
atomic decompositions 
and Banach frames. It allows to transfer certain questions concerning Banach space or interpolation theory 
from the function space to the associated sequence space level, see \cite{FeGr89a, FeGr89b, Ov84}. 
In connection with smoothness spaces of Besov-Lizorkin-Triebel type
the philosophy of this approach is to measure smoothness of a
function in decay properties of the continuous wavelet transform $W_g f$ which is studied in detail in the appendix. 
Indeed, homogeneous Besov and Lizorkin-Triebel type spaces turn out to be 
coorbits of properly chosen spaces $Y$ on the $ax+b$-group $\mathcal{G}$.

There are some more examples according to this abstract theory. One main class of examples refers to the Heisenberg group $\mathbb{H}$, the short-time
Fourier transform and leads to the well-known modulation spaces as coorbits of weighted $L_p(\mathbb{H})$ spaces, see \cite[7.1]{FeGr86} and also
\cite{fegr92-1}.
\subsection{Function spaces on $\mathcal{G}$ }
\label{FSgroup}
Integration on $\mathcal{G}$ will always
be with respect to the left Haar measure $d\mu(x)$. The Haar module on $\mathcal{G}$ is denoted by $\Delta$. 
We define further $L_x F(y) = F(x^{-1}y)$ and
$R_x F(y) = F(yx)$, $x,y \in \mathcal{G}$, the left and right translation operators.
A Banach function space $Y$ on the group $\mathcal{G}$ is supposed to have the
following properties
\begin{description}
  \item(i) $Y$ is continuously embedded in $L_1^{loc}(\mathcal{G})$,
  \item(ii) $Y$ is invariant under left and right translation $L_x$ and $R_x$, which represent in addition
  continuous operators on $Y$,
  \item(iii) $Y$ is solid, i.e., $H\in Y$ and $|F(x)| \leq |H(x)|$ a.e. imply $F\in Y$ and
  $\|F|Y\| \leq \|H|Y\|$.
\end{description}
The continuous weight $w$ is called sub-multiplicative
if $w(xy) \leq w(x) w(y)$ for all $x,y \in \mathcal{G}$. 
The space $L_p^w(\mathcal{G})$, $1\leq p\leq \infty$, of functions $F$ on the group $\mathcal{G}$ is defined via the norm
$$
    \|F|L_p^w(\mathcal{G})\| = \Big(\int\limits_{\mathcal{G}} |F(x)w(x)|^p\,d\mu(x)\Big)^{1/p}\,,
$$
where we use the essential supremum in case $p=\infty$\,. If $w \equiv 1$ then we simply write
$L_p(\mathcal{G})$\,. It is easy to show that these spaces provide left and right translation
invariance if $w$ is sub-multiplicative. Later, in Paragraph \ref{PeFS} we are going to introduce certain
mixed norm spaces where the translation invariance is not longer automatic.

\subsection{Sequence spaces}
\label{sequsp}
\begin{definition} Let $X = \{x_i\}_{i \in I}$ be some discrete set of
points in $\mathcal{G}$ and $V$ be a relatively compact neighborhood
of $e\in \mathcal{G}$\,.
\item(i) $X$ is called $V$-dense if $\mathcal{G} = \bigcup\limits_{i\in I}
x_iV$.
\item(ii) $X$ is called relatively separated if for all compact sets
$K\subset \mathcal{G}$ there exists a constant $C_K$ such that
$$
    \sup\limits_{j\in I} \sharp\{i\in I~:~x_iK \cap x_jK \neq
    \emptyset\} \leq C_K\,.
$$
\item(iii) $X$ is called $V$-well-spread (or simply well-spread) if
it is both relatively separated and $V$-dense for some $V$\,.
\end{definition}
\begin{definition} For a family $X=\{x_i\}_{i\in I}$ which is $V$-well-spread with
respect to a relatively compact neighborhood $V$ of $e\in
\mathcal{G}$ we define the sequence space $Y^b$ and $Y^{\sharp}$ associated to
$Y$ as
\begin{equation}\nonumber
 \begin{split}
    Y^{b} &= \Big\{\{\lambda_i\}_{i\in I}~:~
    \|\{\lambda_i\}_{i\in I}|Y^{b}\| = \Big\|\sum\limits_{i\in I}
    |\lambda_i|\mu(x_iV)^{-1}\chi_{x_iV}|Y\Big\|<\infty\,
    \Big\}\,,\\
    Y^{\sharp} &= \Big\{\{\lambda_i\}_{i\in I}~:~
    \|\{\lambda_i\}_{i\in I}|Y^{\sharp}\| = \Big\|\sum\limits_{i\in I}
    |\lambda_i|\chi_{x_iV}|Y\Big\|<\infty\,\Big\}\,.
  \end{split}
\end{equation}
\end{definition}
\begin{remark} For a well-spread family $X$ the spaces $Y^b$ and $Y^{\sharp}$ do
not depend on the choice of $V$, i.e. different sets $V$ define
equivalent norms on $Y^b$ and $Y^{\sharp}$, respectively\,. For more details on these spaces we refer to 
\cite{FeGr89a}\,.
\end{remark}

\subsection{Coorbit spaces}
Having a Hilbert space $\mathcal{H}$ and an integrable, irreducible, unitary and continuous
representation $\pi:\mathcal{G} \to \mathcal{L}(\mathcal{H})$ then the general voice transform of $f\in \mathcal{H}$ with respect to a fixed atom $g$
is defined as the function $V_gf$ on the group $\mathcal{G}$
given by
\be\label{voice}
    V_gf(x) = \langle \pi(x) g,f\rangle\,,
\ee
where the brackets denote the inner product in $\mathcal{H}$\,.

\begin{definition}\label{H1w}  For a sub-multiplicative weight $w(\cdot)\geq 1$ on $\mathcal{G}$ we define the space
$A_w \subset \mathcal{H}$ of admissible
vectors by
$$
    A_w = \{g\in \mathcal{H}~:~V_g g \in L_1^w(\mathcal{G})\}\,.
$$
If $A_w \neq \{0\}$ and $g\in A_w$ we define further
$$
    \mathcal{H}^1_w(\R) = \{f\in \mathcal{H}~:~\|f|\mathcal{H}^1_w\| = \|V_g f|L^w_1(\mathcal{G})\| <\infty\}\,.
$$
Finally, we denote with $(\mathcal{H}^1_w)^{\sim}$ the canonical anti-dual of
$\mathcal{H}^1_w$, i.e., the space of conjugate linear functionals on $\mathcal{H}^1_w$\,.
\end{definition}
We see immediately that $A_w \subset  \mathcal{H}^1_w \subset \mathcal{H}$. The voice transform \eqref{voice} can now be extended to
$\mathcal{H}_w\times (\mathcal{H}^1_w)^{\sim}$ by the usual dual pairing. The
space $\mathcal{H}^1_w$ can be considered as the space of test functions 
and the reservoir $(\mathcal{H}^1_w)^{\sim}$ as distributions. \\
Let now $Y$ be a space on $\mathcal{G}$ such that (i) - (iii) in Paragraph \ref{FSgroup} hold true.
We define further
\be\label{wy}
    w_Y(x) = \max\{\|L_{x}\|, \|L_{x^{-1}}\|, \|R_{x}\|, \Delta(x^{-1})\|R_{x^{-1}}\| \}\quad,\quad x\in \mathcal{G}\,,
\ee
where the operator norms are considered from $Y$ to $Y$\,.

\begin{definition}\label{defcoob} Let $Y$ be a space on $\mathcal{G}$ satisfying (i)-(iii) in Paragraph \ref{FSgroup}
and let the weight $w(x)$ be given by \eqref{wy}.
  Let further $g\in A_w$. We define the space $\CoY$, which we call coorbit space of $Y$, through
  \be\label{defcoorbit}
      \CoY = \{f\in (\mathcal{H}_w^1)^{\sim}~:~V_g f \in Y\}\quad\mbox{with}\quad
      \|f|\CoY\| = \|V_g f|Y\|\,.
  \ee
\end{definition}
\noindent
The following basic properties are proved for instance in \cite[Thm. 4.5.13]{ra05-6}.
\begin{Theorem}\label{indep}{\em (i)} The space $\CoY$ is a Banach space independent of the analyzing vector $g \in A_w$.\\
      {\em (ii)} The definition of the space $\CoY$ is independent of the reservoir in the
      following sense: Assume that $S \subset \mathcal{H}^1_w$ is a non-trivial locally convex vector
      space which is invariant under $\pi$. Assume further that there exists a non-zero
      vector $g \in S\cap A_w$ for which the reproducing formula
      $$
            V_g f = V_g g \ast V_g f
      $$
      holds true for all $f\in S^{\sim}$. Then we have
      $$
            \CoY = \{f\in (\mathcal{H}_w^1)^{\sim}~:~V_g f \in Y\} = \{f\in S^{\sim}~:~V_g f \in Y\}\,.
      $$
\end{Theorem}

\subsection{Discretizations}
\label{groupdisc}
This section collects briefly the basic facts concerning atomic
(frame) decompositions in coorbit spaces. We are interested in atoms
of type $\{\pi(x_i)g\}_{i\in I}$, where $\{x_i\}_{i\in I} \subset
\mathcal{G}$ represents a discrete subset, whereas $g$ denotes a
fixed admissible analyzing vector.

\begin{definition}\label{atdec}
A family $\{g_i\}_{i\in I}$ in a Banach space $B$ is called an
atomic decomposition for $B$ if there exists a family of bounded
linear functionals $\{\lambda_i\}_{i\in I} \subset B'$ (not necessary
unique) and a Banach sequence space $B^{\sharp} = B^{\sharp}(I)$ such
that:
\begin{description}
    \item(a) We have $\{\lambda_i(f)\}_{i\in I} \in B^{\sharp}$ for all $f\in
    B$ and there exists a constant $C_1>0$ with
    $$
        \|\{\lambda_i(f)\}_{i\in I}|B^{\sharp}\| \leq C_1\|f|B\|\,.
    $$
    \item(b) For all $f\in B$ we have
    $$
        f = \sum\limits_{i\in I} \lambda_i(f)g_i
    $$
    in some suitable topology.
    \item(c) If $\{\lambda_i\}_{i\in I} \in B^{\sharp}$
    then $\sum_{i\in I} \lambda_i g_i \in B$ and there exists a
    constant $C_2>0$ such that
    $$
            \Big\|\sum\limits_{i\in I} \lambda_i g_i|B\Big\| \leq
            C_2 \|\{\lambda_i\}_{i\in I}|B^{\sharp}\|\,.
    $$
\end{description}
\end{definition}

\begin{definition}\label{Bframe} A family $\{h_i\}_{i\in I} \subset B'$ is called a
Banach frame for $B$ if there exists a Banach sequence space
$B^{b} = B^{b}(I)$ and a linear bounded reconstruction
operator $\Theta:B^{b} \to B$ such that:
\begin{description}
    \item(a) We have $\{h_i(f)\}_{i\in I} \in B^{b}$ for all $f\in B$
    and there exist constants $C_1,C_2$ such that
    $$
        C_1\|f|B\| \leq \|\{h_i(f)\}_{i\in I}|B^{b}\| \leq
        C_2\|f|B\|\,,
    $$
    \item(b) and $\Theta(\{h_i(f)\}_{i\in I}) =f$\,.
\end{description}

\end{definition}

\begin{remark} This setting differs slightly from the understanding of
Triebel in \cite{Tr92,Tr06}\,.
\end{remark}
The following abstract result for the atomic decomposition in $\CoY$
is due to Feichtinger and Gr\"ochenig (see \cite[Thm.
6.1]{FeGr89a}).

\begin{Theorem}\label{ADabstr} Let $Y$ be a function space on the group $\mathcal{G}$
satisfying the hypotheses (i)-(iii) from Paragraph \ref{FSgroup} and let $w(x)$ be given by 
\eqref{wy}. Furthermore, the element $g \in A_w$ is supposed to satisfy
\be\label{Vgg}
    \int\limits_{\mathcal{G}} \big(\sup\limits_{y\in xV}
    |\langle\pi(y)g,g\rangle|\big)w(x,t) d\mu(x) < \infty\,.
\ee
Then there exists a neighborhood $U$ of $e\in \mathcal{G}$ and
constants $C_0, C_1>1$ such that for every $U$-well-spread discrete
set $X = \{x_i\}_{i\in I} \subset \mathcal{G}$ the following is
true.
\begin{description}
    \item(i) (Analysis) Every $f\in \CoY$ has a representation
    $$
        f = \sum\limits_{i\in I} \lambda_{i} \pi(x_i)g
    $$
    with coefficients $\{\lambda_i\}_{i\in I}$ depending linearly on
    $f$ and satisfying the estimate
    $$
        \|\{\lambda_{i}\}_{i\in I}|Y^{\sharp}\| \leq C_0\|f|Y\|\,.
    $$
    \item(ii) (Synthesis) Conversely, for any sequence $\{\lambda_i\}_{i\in I} \in Y^{\sharp}$
    the element $f = \sum_{i\in I}\lambda_i \pi(x_i)g$ is in $\CoY$
    and one has
    $$
        \|f|\CoY\| \leq C_1\|\{\lambda_i\}_{i \in I}|Y^{\sharp}\|\,.
    $$
\end{description}
In both cases, convergence takes place in the norm of $\CoY$ if the
finite sequences are norm dense in $Y^{\sharp}$, and in the
weak$^{\ast}$-sense of $(\mathcal{H}_w^1)^{\sim}$ otherwise.
\end{Theorem}
\begin{remark} According to Definition \ref{atdec} the family $\{\pi(x_i)g\}_{i\in I}$
  represents an atomic decomposition for $\CoY$.
\end{remark}
\begin{Theorem}\label{sabound} Under the same assumptions as in Theorem \ref{ADabstr} the system
$\{\pi(x_i) g\}_{i\in I}$ represents a Banach frame for $\CoY$, i.e.,
  $$
        \|f|\CoY\| \asymp \|\{\langle \pi(x_i)g, f  \rangle\}|Y^b\|\quad,\quad f\in \CoY\,.
  $$
\end{Theorem}
\noindent
The following powerful result goes back to Gr\"ochenig \cite{Gr88}
and was generalized by Rauhut \cite{ra05-3}.
\begin{Theorem}\label{wbases} Suppose that the functions $g_r, \gamma_r$, $r=1,...,n$, satisfy
\eqref{Vgg}. Let $X = \{x_i\}_{i\in I}$ be a well-spread set such
that
\be\label{waveletexp}
    f = \sum\limits_{r=1}^n \sum\limits_{i\in I} \langle
    \pi(x_i)\gamma_r, f\rangle \pi(x_i)g_r
\ee
for all $f\in \mathcal{H}$\,. Then expansion \eqref{waveletexp} extends to all
$f\in \CoY$. Moreover, $f\in (\mathcal{H}_w^1)^{\sim}$ belongs to $\CoY$ if and only
if $\{\langle \pi(x_i)\gamma_r,f\rangle\}_{i\in I}$ belongs to
$Y^{b}$ for each $r=1,...,n$\,. The convergence is considered
in $\CoY$ if the finite sequences are dense in $Y^{b}$. In
general we have weak$^{\ast}$-convergence.
\end{Theorem}

\bproof The proof of this result relies on the fact, that there
exists an atomic decomposition $\{\pi(y_i)g\}_{i\in I}$ by Theorem
\ref{ADabstr} with a certain $g$ satisfying \eqref{Vgg} and a corresponding sequence of points 
$Z = \{y_i\}_{i\in I}$\,. This has
to be combined with Theorem \ref{sabound} and Theorem \ref{ADabstr}/(ii) and we are done. See
\cite{Gr88} for the details.
\eproof.

\section{Coorbit spaces on the $ax+b$-group }
\label{ax+b}

Let $\mathcal{G} = \mathbb{R}^d \rtimes \mathbb{R}_+^*$ the $d$-dimensional $ax+b$-group. Its multiplication
is given by
$$
  (x,t)(y,s) = (x+ty,st)\,.
$$
The left Haar measure $\mu$ on $\mathcal{G}$ is given by
$d\mu(x,t) = dx\,dt/t^{d+1}$, the Haar module is $\Delta(x,t) = t^{-d}$.
Giving a function $F$ on $\mathcal{G}$ the left and right
translation $L_y = L_{(y,r)}$ and $R_y = R_{(y,r)}$ are given by
$$
    L_{(y,r)}F(x,t) = F((y,r)^{-1}(x,t)) = F\Big(\frac{x-y}{r},\frac{t}{r}\Big)
$$
and
$$
  R_{(y,r)}F(x,t) = F((x,t)(y,r)) = F(x+ty,rt)\,.
$$

\subsection{Peetre type spaces on $\mathcal{G}$}
\label{PeFS}
The present paragraph is devoted to the definition of certain mixed norm spaces on the group.
Such spaces have been considered in various papers, see
\cite{CoMeSt85, FeGr86, Gr88, Gr91}. In particular, so-called tent spaces
have some important applications in harmonic analysis. Indeed, it is possible 
to recover Lizorkin-Triebel spaces as coorbits of tent spaces.\\
Here we use a different approach and define a new scale of function spaces on the group $\mathcal{G}$. We call
them Peetre type spaces since a quantity related to the Peetre maximal function \eqref{Peetrefunct} is involved in its definition. It turned
out that they are straight forward to handle in connection with translation invariance. In contrast to the tent space approach
they represent the more natural choice for considering Lizorkin-Triebel spaces as coorbits. Additionally,
they seem to be suitable for inhomogeneous spaces and
more general situations like weighted spaces and general $2$-microlocal spaces, which will be studied in a further contribution to the subject.

\begin{definition}\label{sponG}
  Let $s\in \re$, $0<p,q\leq \infty$, and $a>0$. We define the spaces $\dot{L}_{p,q}^s(\mathcal{G})$, 
  $\dot{T}^{s}_{p,q}(\mathcal{G})$, and $\dot{P}_{p,q}^{s,a}(\mathcal{G})$ on the group $\mathcal{G}$ via the finiteness of the following (quasi-)norms 
  \begin{equation}
  \label{tent}
  \begin{split}
     \|F|\dot{L}_{p,q}^s(\mathcal{G})\| &= \Big(\int\limits_{0}^{\infty}t^{-sq}\|F(\cdot,t)|L_p(\R)\|^q \frac{dt}{t^{d+1}}\Big)^{1/q}\,,\\
     \|F|\dot{T}^{s}_{p,q}(\mathcal{G})\| &= \Big\|\Big(\int\limits_{0}^{\infty} t^{-sq}\int\limits_{B(0,t)}|F(x+z,t)|^q\,dz\frac{dt}{t^{d+1}}\Big)^{1/q}|L_p(\R)\Big\|\,,\\
     \|F|\dot{P}_{p,q}^{s,a}(\mathcal{G})\| &= \Big\|\Big(\int\limits_{0}^{\infty} t^{-sq}\Big[\sup\limits_{y\in \R} \frac{|F(x+y,t)|}{(1+|y|/t)^a}\Big]^q\frac{dt}{t^{d+1}}\Big)^{1/q}|L_p(\R)\Big\|\,,
  \end{split}
  \end{equation} 
    using the usual modification in case $q=\infty$.
\end{definition}

\begin{proposition}\label{rti} The spaces $\dot{L}^s_{p,q}(\mathcal{G}),\dot{T}^s_{p,q}(\mathcal{G})$ and $\dot{P}^s_{p,q}(\mathcal{G})$ are
left and right translation invariant. Precisely, we have
\be\nonumber
  \begin{split}
    \|L_{(z,r)}:\dotspace{L}{s} \to \dotspace{L}{s}\| &= r^{d(1/p-1/q)-s}\,,\\
    \|R_{(z,r)}:\dotspace{L}{s} \to \dotspace{L}{s}\| &= r^{s+d/q}\,,\\
    \|L_{(z,r)}:\dotspace{T}{s} \to \dotspace{T}{s}\| &= r^{d/p-s}\,,\\
    \|R_{(z,r)}:\dotspace{T}{s} \to \dotspace{T}{s}\| &\leq Cr^{d/q+s}\max\{1,r^{-b}(1+|z|)^{b}\}\,,
  \end{split}
\ee
where $b>0$ is a constant depending on $d$, $p$ and $q$. Furthermore, we have
\be\nonumber
    \begin{split}
    \|L_{(z,r)}:\dotspace{P}{s,a} \to \dotspace{P}{s,a}\| &= r^{d(1/p-1/q)-s}\,,\\
    \|R_{(z,r)}:\dotspace{P}{s,a} \to \dotspace{P}{s,a}\| &\leq r^{s+d/q}\max\{1,r^{-a}\}(1+|z|)^a\,.
    \end{split}
\ee
\end{proposition}
\bproof {\em Step 1.} The left and right translation invariance of
$\dotspace{L}{s}$ and $\dotspace{T}{s}$ was shown in \cite[Lem. 4.7.10]{ra05-6}.\frei\\\newline
{\em Step 2.} Let us consider $\dotspace{P}{s,a}$. Clearly, we have for $F\in \dotspace{P}{s,a}$
\be\nonumber
    \begin{split}
      \|L_{(z,r)}F|\dotspace{P}{s}\| &= \Big\|\Big(\int\limits_{0}^{\infty}t^{-sq}\Big[
      \sup\limits_{y\in\R}\frac{|F((x+y-z)/r,t/r)|}{(1+|y|/t)^a}\Big]^q\frac{dt}{t^{d+1}}\Big)^{1/q}\Big|L_p(\R)\Big\|\\
      &= r^{d/p}\Big\|\Big(\int\limits_{0}^{\infty}t^{-sq}\Big[
      \sup\limits_{y\in\R}\frac{|F(x+y,t/r)|}{(1+r|y|/t)^a}\Big]^q\frac{dt}{t^{d+1}}\Big)^{1/q}\Big|L_p(\R)\Big\|\\
      &=r^{d(1/p-1/q)-s}\Big\|\Big(\int\limits_{0}^{\infty}t^{-sq}\Big[
      \sup\limits_{y\in\R}\frac{|F(x+y,t)|}{(1+|y|/t)^a}\Big]^q\frac{dt}{t^{d+1}}\Big)^{1/q}\Big|L_p(\R)\Big\|\,.
    \end{split}
\ee
Hence, we obtain
$$
    \|L_{(z,r)}:\dotspace{P}{s,a} \to \dotspace{P}{s}\| = r^{d(1/p-1/q)-s}\,.
$$
The right translation invariance is obtained by \be\nonumber
    \begin{split}
      \|R_{(z,r)}F|\dotspace{P}{s,a}\| &= \Big\|\Big(\int\limits_{0}^{\infty}t^{-sq}\Big[
      \sup\limits_{y\in\R}\frac{|F(x+tz+y,tr)|}{(1+|y|/t)^a}\Big]^q\frac{dt}{t^{d+1}}\Big)^{1/q}\Big|L_p(\R)\Big\|\\
      &= \Big\|\Big(\int\limits_{0}^{\infty}t^{-sq}\Big[
      \sup\limits_{y\in\R}\frac{|F(x+y,tr)|}{(1+|y-tz|/t)^a}\Big]^q\frac{dt}{t^{d+1}}\Big)^{1/q}\Big|L_p(\R)\Big\|\\
      &= r^{s+d/q}\Big\|\Big(\int\limits_{0}^{\infty}t^{-sq}\Big[
      \sup\limits_{y\in\R}\frac{|F(x+y,t)|}{(1+|y-tz|r/t)^a}\Big]^q\frac{dt}{t^{d+1}}\Big)^{1/q}\Big|L_p(\R)\Big\|.
    \end{split}
\ee
Observe that
\be\nonumber
    \sup\limits_{y\in\R}\frac{|F(x+y,t)|}{(1+|y-tz|r/t)^a} =
    \sup\limits_{y\in\R}\Big[\frac{|F(x+y,t)|}{(1+|y|/t)^a}\cdot \frac{(1+|y|/t)^a}{(1+|y-tz|r/t)^a}\Big]
\ee
and
$$
  \frac{(1+|y|/t)^a}{(1+|y-tz|r/t)^a}\leq \frac{(1+|y-tz|/t+|z|)^a}{(1+|y-tz|r/t)^a} =
  \frac{(1+|y-tz|/t)^a)(1+|z|)^a}{(1+|y-tz|r/t)^a} \,.
$$
This yields
$$
  \sup\limits_{y\in\R}\frac{|F(x+y,t)|}{(1+|y-tz|r/t)^a} \leq \max\{1,r^{-a}\}(1+|z|)^a \sup\limits_{y\in\R}\frac{|F(x+y,t)|}{(1+|y|/t)^a}\,
$$
and consequently
$$
    \|R_{(z,r)}:\dotspace{P}{s,a} \to \dotspace{P}{s,a}\| \leq r^{s+d/q}\max\{1,r^{-a}\}(1+|z|)^a\,.
$$
\eproof
\begin{remark}
Note, that we did neither use the translation invariance of the Lebesgue
measure nor any change of variable in order to prove the right translation invariance of
$\dotspace{P}{s,a}$. This gives room for further generalizations, i.e., replacing
the space $L_p(\R)$ by some weighted Lebesgue space $L_p(\R,\omega)$ for instance.
\end{remark}

\subsection{New old coorbit spaces}
We start with $\mathcal{H}=L_2(\R)$ and the representation
$$
    \pi(x,t) = T_x \mathcal{D}_t^{L_2}\,,
$$
where $T_xf = f(\cdot-x)$ and $\mathcal{D}_t^{L_2}f = t^{-d/2}f(\cdot/t)$ has been already defined in Paragraph \ref{sectmax}.
This representation is unitary, continuous and square integrable on $\mathcal{H}$ but not irreducible. However, if we
restrict to radial functions $g\in L_2(\R)$ then $\mbox{span} \{\pi(x,t)g~:(x,t) \in \mathcal{G}\}$ is dense
in $L_2(\R)$. Another possibility to overcome this obstacle is to extend the group by $SO(d)$, which is more or less equivalent,
see \cite{FeGr86, FeGr89a} for details. The voice transform in this special situation is represented by the so-called continuous 
wavelet transform $W_g f$ which we study in detail in Paragraph \ref{SectCWT} in the appendix. 

Recall the abstract definition of the space $\mathcal{H}^1_w$ and $\mathcal{A}_w$ from Definition \ref{H1w}.
The following result implied by our
Lemma \ref{help1} on the decay of the continuous wavelet transform. It states under which conditions on the weight $w$
the space $\mathcal{H}^1_w$ is nontrivial.
\begin{lemma}\label{emb} If the weight function $w(x,t)\geq 1$ satisfies the condition
    $$
        w(x,t) \leq (1+|x|)^r(t^s+t^{-s'})
    $$
    for some $r,s,s'\geq 0$ then
    $$
          \mathcal{S}_0(\R) \hookrightarrow \mathcal{H}^1_w\,.
    $$
\end{lemma}
This is a kind of minimal condition which is needed in order to define
coorbit spaces in a reasonable way. Instead of $(H^1_w)^{\sim}$ one may use
$\mathcal{S}_0'(\R)$ as reservoir and a radial $g\in \mathcal{S}_0(\R)$ as analyzing vector. Considering
\eqref{wy} we have to restrict to such function spaces $Y$ on $\mathcal{G}$ satisfying (i),(ii),(iii) in
Paragraph \ref{FSgroup} where additionally 
\begin{description}
 \item(iv)
 $$
    w(x,t) = w_Y(x,t) \lesssim (1+|x|)^r(t^s+t^{-s'})
 $$
\end{description}
holds true for some $r,s,s'\geq 0$\,. The following theorem shows, how the spaces of
Besov-Lizorkin-Triebel type from Section \ref{clFS} can be recovered as coorbit spaces with respect to $\mathcal{G}$.
\begin{Theorem}\label{coohom}
  \begin{description}
      \item(i) For $1\leq p,q\leq \infty$ and $s\in \re$ we have
      $$
          \dot{B}^{s}_{p,q}(\R) = \mbox{Co}\dot{L}^{s+d/2-d/q}_{p,q}(\mathcal{G})\,,
      $$
      \item(ii) for $1\leq p<\infty$, $1\leq q\leq \infty$ and $s\in \re$ we have
      $$
          \dot{F}^{s}_{p,q}(\R) = \mbox{Co}\dot{T}^{s+d/2}_{p,q}(\mathcal{G})
      $$
      \item(iii) and if additionally $a>\frac{d}{\min\{p,q\}}$ we obtain
      $$
          \dot{F}^{s}_{p,q}(\R) = \mbox{Co}\dot{P}^{s+d/2-d/q,a}_{p,q}(\mathcal{G})\,.
      $$
  \end{description}
\end{Theorem}
\bproof Theorem \ref{coohom} is a
direct consequence of Definition \ref{defcoob}, formula \eqref{CWTconv}, Proposition \ref{rti}, Theorems
\ref{main1b}, \ref{main1bb} and the abstract result in Theorem \ref{indep}. \eproof

\begin{remark} {\em (a)} The assertions (i) and (ii) are not new. They appear for instance in \cite{FeGr86, Gr88, Gr91} and 
  rely on the characterizations given by Triebel in \cite{Tr88} and \cite[2.4, 2.5]{Tr92}, see in particular 
  \cite[2.4.5]{Tr92} for the variant in terms of tent spaces which were invented in \cite{CoMeSt85}. From the deep result
  in \cite[Prop. 4]{CoMeSt85} it follows that $\dot{T}^{s}_{p,q}(\mathcal{G})$ are translation invariant Banach function spaces on $\mathcal{G}$, 
  which makes them feasible for coorbit space theory\\
  
  {\em (b)} Assertion (iii) is indeed new and makes the rather complicated tent spaces $\dot{T}^{s}_{p,q}(\mathcal{G})$
  obsolete for this issue. We showed that $Y = P^{s,a}_{p,q}(\mathcal{G})$ is a much better choice
  since the right translation invariance is immediate and gives more transparent estimates for its norm. Once we are interested
  in reasonable conditions for atomic decompositions this is getting important, see Section \ref{Wframes}.
\end{remark}

\subsection{Sequence spaces}
\label{sequhom}
In the sequel we consider a compact neighborhood
of the identity element in $\mathcal{G}$ given by $\mathcal{U} =
[0,\alpha]^d \times [\beta^{-1},1]$, where $\alpha>0$ and $1<\beta$.
Furthermore, we consider the discrete set of points
$$
  \{x_{j,k} = (\alpha k \beta^{-j},\beta^{-j})~:~j\in \zz, k\in \Z\}\,.
$$
This family is $\mathcal{U}$-well-spread. Indeed,
$$
   x_{j,k}\mathcal{U} = Q_{j,k} \times [\beta^{-(j+1)},\beta^{-j}]\,,
$$
where
$$
    Q_{j,k} = [\alpha k_1\beta^{-j},\alpha (k_1+1)\beta^{-j}]\times \cdots \times[\alpha k_d\beta^{-j},\alpha (k_d+1)\beta^{-j}]\,.
$$
Note that in this case the spaces $Y^{\sharp}$ and $Y^b$ coincide.
We will further use the notation
$$
    \chi_{j,k}(x) = \left\{\begin{array}{rcl}
                           1 &:& x\in Q_{j,k}\\
                           0 &:& \mbox{otherwise}
                        \end{array}\right..
$$
\begin{definition} Let $Y$ be a function space on $\mathcal{G}$ as above. We put
   $$
      Y^{\sharp}(\alpha,\beta) = \{\{\lambda_{j,k}\}_{j,k}~:~\|\{\lambda_{j,k}\}_{j,k}|Y^{\sharp}(\alpha,\beta)\|<\infty\}\,,
   $$
   where
   $$
      \|\{\lambda_{j,k}\}_{j,k}|Y^{\sharp}(\alpha,\beta)\| = \Big\|\sum\limits_{j,k}|\lambda_{j,k}|\chi_{j,k}(x)\chi_{[\beta^j,\beta^ {j+1}]}(t)|Y\Big\|\,.
   $$
\end{definition}

\begin{Theorem} Let $1\leq p,q \leq \infty$, $s\in \re$ and $a > d/\min\{p,q\}$. Then
$$
    \|\{\lambda_{j,k}\}_{j,k}|(\dot{P}^{s,a}_{p,q})^{\sharp}(\alpha,\beta)\| \asymp
    \Big\|\Big(\sum\limits_{\ell\in \zz}\sum\limits_{k \in \Z}\beta^{\ell (s+d/q)q }|\lambda_{\ell,k}|^q\chi_{\ell,k}(x)\,
       \Big)^{1/q}|L_{p}(\R)\Big\|
$$
and
$$
    \|\{\lambda_{j,k}\}_{j,k}|(\dot{L}^{s}_{p,q})^{\sharp}(\alpha,\beta)\| \asymp
    \Big(\sum\limits_{\ell\in \zz}\beta^{\ell (s+d/q-d/p)q }\Big(\sum\limits_{k \in \Z}|\lambda_{\ell,k}|^p\Big)^{q/p}\,
       \Big)^{1/q}\,.
$$

\end{Theorem}

\bproof We prove the first statement. The proof for the second one is even simpler.
Let
$$
    F(x,t) = \sum\limits_{j,k}|\lambda_{j,k}|\chi_{j,k}(x)\cdot\chi_{[\beta^{-(j+1)},\beta^{-j}]}(t)\,.
$$
Discretizing the integral over $t$ by $t \asymp \beta^{-\ell}$ we
obtain \be\label{eq-4}
    \begin{split}
       \|F|Y\| &= \Big\|\Big(\int\limits_{0}^{\infty}t^{-sq}\Big[\sup\limits_{w}\frac{|F(x+w,t)|}{(1+|w|/t)^a}\Big]^q
       \frac{dt}{t^{d+1}}\Big)^{1/q}|L_p(\R)\Big\|\\
       &\asymp \Big\|\Big(\sum\limits_{\ell\in \zz}\beta^{\ell (s+d/q) q}\int\limits_{\beta^{-(\ell+1)}}^{\beta^{-\ell}}\Big[\sup\limits_{w}\frac{|F(x+w,t)|}{(1+\beta^{\ell}|w|)^a}\Big]^q\,\frac{dt}{t}\Big)^{1/q}|L_p(\R)\Big\|\,.
    \end{split}
\ee
With $t\in [\beta^{-(\ell+1)},\beta^{-\ell}]$ we observe
\be\label{eq-5}
    F(x,t) = \sum\limits_k|\lambda_{\ell,k}|\chi_{\ell,  k}(x)
\ee
and estimate
\be\label{eq-1}
    \begin{split}
       \|F|Y\| \leq  \Big\|\Big(\sum\limits_{\ell\in \zz}\beta^{\ell (s+d/q) q}\sup\limits_{w\in \R}\frac{1}{(1+\beta^{\ell}|w|)^a}\Big[\sum\limits_{k}|\lambda_{\ell,k}|\chi_{\ell,k}(x+w)\Big]^q\,
       \Big)^{1/q}|L_p(\R)\Big\|\,.
    \end{split}
\ee In order to include also the situation $\min\{p,q\}\leq 1$ we use
the following trick. Obviously, we can rewrite and estimate
\eqref{eq-1} with $0<r<1$ in the following way \be\label{eq-2}
    \begin{split}
       \|F|Y\| \leq  \Big\|\Big(\sum\limits_{\ell\in \zz}\Big[
       \sum\limits_{k}\beta^{\ell (s+d/q)r}|\lambda_{\ell,k}|^r\sup\limits_{w\in \R}\frac{\chi_{\ell,k}(x+w)}{(1+\beta^{\ell}|w|)^{ar}}\Big]^{q/r}\,
       \Big)^{1/q}|L_p(\R)\Big\|\,.
    \end{split}
\ee
We continue with the useful estimate
\be\label{eq-0}
      \sup\limits_{w}\frac{|\chi_{\ell,k}(x+w)|}{(1+\beta^{\ell}|w|)^{ar}} \lesssim
      \frac{1}{(1+\beta^{\ell}|x-k\beta^{-\ell}|)^{ar}}\lesssim \Big(\chi_{\ell,k}(\cdot) \ast \frac{\beta^{\ell d}}{(1+\beta^{\ell}|\cdot|)^{ar}}\Big)(x).
\ee
Indeed, the first estimate is obvious. Let us establish the second one
\be\label{eq-34}
  \begin{split}
     \Big(\chi_{\ell,k}(\cdot) \ast \frac{1}{(1+\beta^{\ell}|\cdot|)^{ar}}\Big)(x) &=
     \int\limits_{\substack{|y_i-k_i\beta^{-\ell}|\leq \beta^{-\ell}\\i=1,...,d}}
     \frac{1}{(1+\beta^{\ell}|x-y|)^{ar}} \,dy\\
     &\gtrsim \int\limits_{|y|\leq c\beta^{-\ell}}
    \frac{1}{(1+\beta^{\ell}|x-k\beta^{-\ell}-y|)^{ar}}\,dy\\
    &\gtrsim \int\limits_{|y|\leq c\beta^{-\ell}}
    \frac{1}{(1+\beta^{\ell}|x-k\beta^{-\ell}|+\beta^{\ell}|y|)^{ar}}\,dy\\
    &\gtrsim \beta^{-\ell d}\int\limits_0^{1}
    \frac{u^{d-1}}{(1+\beta^{\ell}|x-k\beta^{-\ell}|+u)^{ar}}\,du\\
    &\gtrsim \frac{\beta^{-\ell d}}{(1+\beta^{\ell}|x-k\beta^{-\ell}|)^{ar}}\,.
   \end{split}
\ee
Note, that the functions
$$
    g_{\ell}(x)  = \frac{\beta^{\ell d}}{(1+\beta^{\ell}|\cdot|)^{ar}}\quad
$$
belong to $L_1(\R)$ with uniformly bounded norm,
where we need that $ar>d$\,. Putting \eqref{eq-34} and \eqref{eq-0} into \eqref{eq-2} we obtain
\be\nonumber
    \begin{split}
       \|F|Y\|^r \leq  \Big\|\Big(\sum\limits_{\ell\in \zz}\Big[g_{\ell}\ast\sum\limits_{k \in \Z}\beta^{\ell (s+d/q)r }|
       \lambda_{\ell,k}|^r\chi_{\ell,k})(x)\Big]^{q/r}\,
       \Big)^{r/q}|L_{p/r}(\R)\Big\|\,.
    \end{split}
\ee
Now we are in a position to use the majorant property of the Hardy-Littlewood maximal operator
(see Paragraph \ref{sectmax} and \cite[Chapt. 2]{StWe71}), which states that a convolution of a
function $f$ with a $L_1(\R)$-function (having norm one) can be estimated from above by the Hardy-Littlewood
maximal function of $f$. We choose $r<\min\{p,q\}$ and apply Theorem \ref{feffstein} for the $L_{p/r}(\ell_{q/r})$ situation. This
gives
\be\nonumber
    \begin{split}
       \|F|Y\|^r &\lesssim \Big\|\Big(\sum\limits_{\ell\in \zz}\Big[\sum\limits_{k \in \Z}\beta^{\ell (s+d/q)r }|
       \lambda_{\ell,k}|^r\chi_{\ell,k}(x)\Big]^{q/r}\,
       \Big)^{r/q}|L_{p/r}(\R)\Big\|\\
       &=  \Big\|\Big(\sum\limits_{\ell\in \zz}\sum\limits_{k \in \Z}\beta^{\ell (s+d/q)q }|\lambda_{\ell,k}|^q\chi_{\ell,k}(x)\,
       \Big)^{1/q}|L_{p}(\R)\Big\|^r\,
    \end{split}
\ee
and finishes the upper estimate. Both conditions, $ar>d$ and $r<\min\{p,q\}$, are compatible if $a>d/\min\{p,q\}$ is assumed at the beginning. \\
For the estimate from below we go back to \eqref{eq-4}
and observe
$$
  \sup\limits_{w}\frac{|F(x+w,t)|}{(1+\beta^{\ell}|w|)^a} \geq |F(x,t)|\,,
$$
which results in
\be\nonumber
  \begin{split}
    \|F|Y\| &\gtrsim \Big\|\Big(\sum\limits_{\ell\in \zz}\beta^{\ell (s+d/q) q}\int\limits_{\beta^{-(\ell+1)}}^{\beta^{-\ell}}|F(x,t)|^q\,\frac{dt}{t}\Big)^{1/q}|L_p(\R)\Big\|\,.\\
    &
  \end{split}
\ee
A further use of \eqref{eq-5} gives finally
$$
  \|F|Y\| \gtrsim \Big\|\Big(\sum\limits_{\ell\in \zz}\sum\limits_{k \in \Z}\beta^{\ell (s+d/q)q }|\lambda_{\ell,k}|^q\chi_{\ell,k}(x)\,
   \Big)^{1/q}|L_{p}(\R)\Big\|\,.
$$
The proof is complete.\eproof

\subsection{Atomic decompositions}
The following theorem is a direct consequence of the abstract results
in Theorems \ref{ADabstr}, \ref{sabound}.

\begin{Theorem} Let $1\leq p,q\leq \infty$, $a>d/\min\{p,q\}$ and $s\in \re$. Let further
$g\in \mathcal{S}_0(\R)$ be a radial function. Then there exist numbers $\alpha_0>0$ and $\beta_0>1$
such that for all $0<\alpha\leq \alpha_0$ and $1< \beta\leq \beta_0$
the family
$$
   \{g_{j,k}\}_{j\in\zz,k\in\Z} = \{T_{\alpha k \beta^j}\mathcal{D}^{L_2}_{\beta^j}g\}
$$
has the following properties:
\begin{description}
   \item(i) $\{g_{j,k}\}_{j\in\zz,k\in\Z}$ forms a Banach frame for
   $\mbox{Co}\dotspace{L}{s}$ and $\mbox{Co}\dotspace{P}{s,a}$, i.e., we have a dual frame
   $\{e_{j,k}\}_{j\in\zz,k\in \Z}\subset \mathcal{S}_0(\R)$ with $f=\sum_{j\in\zz,k\in\Z} \langle g_{j,k},f \rangle
   e_{j,k}$ and the norm equivalences
   $$
        \|f|\mbox{Co}\dotspace{L}{s}\| \asymp \|\langle
        g_{j,k},f\rangle|(\dotspace{L}{s})^{\sharp}(\alpha,\beta)\|\quad,\quad
        f\in \mbox{Co}\dotspace{L}{s}
   $$
   as well as
   $$
        \|f|\mbox{Co}\dotspace{P}{s,a}\| \asymp \|\langle
        g_{j,k},f\rangle|(\dotspace{P}{s,a})^{\sharp}(\alpha,\beta)\|\quad,\quad
        f\in \mbox{Co}\dotspace{P}{s,a}\,.
   $$

   \item(ii) $\{g_{j,k}\}_{j\in\zz,k\in\Z}$ is an atomic
   decomposition, i.e.,
   for $f\in \mbox{Co}\dotspace{P}{s,a}$ we have a (not necessary unique) decomposition
   $\sum_{j\in\zz,k\in \Z} \lambda_{j,k}(f)g_{j,k}$ such that
   $$
        \|\{\lambda_{j,k}(f)\}_{j,k}|(\dotspace{P}{s,a})^{\sharp}(\alpha,\beta)\|
        \lesssim \|f|\mbox{Co}\dotspace{P}{s,a}\|\,.
   $$
   Conversely, if $\{\lambda_{j,k}\}_{j\in \zz,k\in \Z}
   \in (\dotspace{P}{s,a})^{\sharp}(\alpha,\beta)$ then
   $f = \sum_{j\in \zz,k\in \Z} \lambda_{j,k} g_{j,k}$ converges
   and belongs to $\mbox{Co}\dotspace{P}{s,a}$ and moreover,
   $$
        \|f|\mbox{Co}\dotspace{P}{s,a}\| \lesssim \|\{\lambda_{j,k}\}_{j,k}|
        (\dotspace{P}{s,a})^{\sharp}(\alpha,\beta)\|\,
   $$
   (analogously for $\mbox{Co}\dotspace{L}{s}$). Convergence is considered in the strong topology if the
   finite sequences are dense in $(\dotspace{P}{s,a})^{\sharp}(\alpha,\beta)$ and in the weak$^{\ast}$-topology
   otherwise.
\end{description}
\end{Theorem}

\begin{remark} {\em (i)} Since the analyzing function or atom $g$ can be chosen
arbitrarily we allow more flexibility here than in the results
given in Frazier/Jawerth \cite{FrJa90} and Triebel \cite{Tr92, Tr06}.\\
{\em (ii)} Instead of regular families of sampling points
$(\alpha\beta^{-j}k,\beta^{-j})$ rather irregular families of points in
$\mathcal{G}$ are allowed as long as they are distributed
sufficiently dense, see Theorem \ref{ADabstr}.
\end{remark}

\subsection{Wavelet frames}
\label{Wframes}
In the sequel we consider wavelet bases on $\R$ in the sense of Lemma \ref{dwavelet} in the appendix. We have given
an orthonormal scaling function $\Psi^0$ and the associated wavelet $\Psi^1$ on $\re$ and consider the tensor products $\Psi^c$, $c\in E$. Our aim is to 
specify, i.e., give sufficient conditions to $\Psi^0$, $\Psi^1$, such that \eqref{wavelet} represents an unconditional basis
in $B^s_{p,q}(\R)$ and $F^s_{p,q}(\R)$, respectively. We intend to apply our abstract Theorem \ref{wbases} and need therefore
to have \eqref{Vgg} for all functions $\Psi^c$. To ensure this we impose certain smoothness $(S_K)$, decay $(D)$, and
moment conditions $(M_L)$ to $\Psi^1$ and $\Psi^0$, which are specified in Definition \ref{basedef}.

\begin{proposition}\label{propwiener} Let $L\in \N$, $K>0$, and $\Psi^0$ be an orthogonal scaling function with associated wavelet $\Psi^1$ on $\re$.
The function $\Psi^0$ is supposed to satisfy $(D)$ and $(S_K)$ and $\Psi^1$ is supposed to satisfy $(D)$, $(S_K)$ and $(M_{L-1})$\,. Let 
$V = [-1,1]^d \times (1/2,1] \subset \mathcal{G}$ a neighborhood of the identity $e\in \mathcal{G}$. 
Suppose further that for $r_1,r_2\in \re$ the weight $w(x,t)$ is given by
$$
    w(x,t) = (1+|x|)^{v}(t^{r_2}+t^{-r_1})\quad,\quad (x,t) \in \mathcal{G}\,.
$$
If now
\be\label{f19}
      r_1 < \min\{L,K\}-d/2\quad,\quad r_2 < \min\{L,K\}+d/2-v
\ee
then we have 
$$
    \int\limits_{\R}\int\limits_{0}^{\infty} \sup\limits_{(y,s)\in (x,t)V}|\langle \pi(y,s)\Psi^c,\Psi^c\rangle|
    w(x,t) \frac{dt}{t^{d+1}}\,dx<\infty\,.
$$
\end{proposition}
\bproof With Lemma \ref{help1} we obtain for $W_{\Psi^1}\Psi^1$ the following
estimates
\be\nonumber
    |(W_{\Psi^1}\Psi^1)(s,t)| \lesssim
    \frac{t^{\min\{L,K\}+1/2}}{(1+t)^{2\min\{L,K\}+1}}\cdot \frac{1}{(1+|s|/(1+t))^{N}}
\ee
And in addition
\be\nonumber
    |(W_{\Psi^i}\Psi^i)(s,t)| \lesssim
    \frac{t^{1/2}}{(t+1)}\frac{1}{(1+|s|/(1+t))^{N}}\quad,\quad i=1,2.
\ee
Hence, for any $c\in E$ the tensor product structure gives (assume without restriction that $c_d=1$)
\be\nonumber
  \begin{split}
   &|W_{\Psi^c}\Psi^c(x,t)| \lesssim
   \frac{t^{\min\{L,K\}}}{(1+t)^{2\min\{L,K\}}}\cdot \frac{t^{d/2}}{(1+t)^d}
   \prod\limits_{i=1}^{d}\frac{1}{(1+|x_i|/(1+t))^{N}}\,.
  \end{split}
\ee
The expression $\sup\limits_{(y,s) \in (x,t)V} |W_{\Psi^c}\Psi^c(y,s)|$ can be estimated similar
\be\nonumber
  \begin{split}
   \sup\limits_{(y,s) \in (x,t)V} |W_{\Psi^c}\Psi^c(y,s)|
   &=\sup\limits_{\substack{|y_i-x_i|\leq t \\ t/2\leq s\leq t}} |W_{\Psi^c}\Psi^c(y,s)|\\
   &\lesssim \frac{t^{\min\{L,K\}}}{(1+t)^{2\min\{L,K\}}}\frac{t^{d/2}}{(1+t)^d}
   \prod\limits_{i=1}^{d}\frac{1}{(1+|x_i|/(1+t))^{N}}\\
   &\lesssim \frac{t^{\min\{L,K\}}}{(1+t)^{2\min\{L,K\}}}\frac{t^{d/2}}{(1+t)^d}
   \frac{1}{(1+|x|/(1+t))^{N}}
  \end{split}
\ee
Fubini's theorem and a change of variable yields
\be\nonumber
    \int\limits_{\R}\int\limits_{0}^{\infty} \sup\limits_{(y,s)\in (x,t)V}|\langle \pi(y,s)\Psi^c,\Psi^c\rangle|
    w(x,t) \frac{dt}{t^{d+1}}\,dx
    \lesssim \int\limits_{0}^{\infty} \frac{t^{\min\{L,K\}}}{(1+t)^{2\min\{L,K\}}}t^{d/2}
    \cdot(1+t)^v(t^{r_2}+t^{-r_1})\frac{dt}{t^{d+1}}\,.
\ee
Finally it is easy to see that the latter is finite if the conditions in \eqref{f19} are valid.
This proves Proposition \ref{propwiener}.\eproof

\begin{Theorem}\label{wavdec} Let $L\in \N$, $K>0$, and $\Psi^0$ be an orthogonal scaling function with associated wavelet $\Psi^1$ on $\re$.
The function $\Psi^0$ is supposed to satisfy $(D)$ and $(S_K)$ and $\Psi^1$ is supposed to satisfy $(D)$, $(S_K)$ and $(M_{L-1})$\,. 
\begin{description}
    \item(a) If $1\leq p,q \leq \infty$ and
    $$
            -\min\{L,K\}+\frac{d}{p}<s<\min\{L,K\}-d\Big(1-\frac{1}{p}\Big)
    $$
    then \eqref{wavelet} is a Banach frame for $\dot{B}^s_{p,q}(\R)$ in the sense
    of \eqref{waveletexp}.
    \item(b) If $1\leq p <\infty$, $1 \leq q\leq \infty$, and
    $$
            -\min\{L,K\}+2d\max\Big\{\frac{1}{p},\frac{1}{q}\Big\}<s<\min\{L,K\}-d\max\Big\{\frac{1}{p},\frac{1}{q},1-\frac{1}{p}\Big\}
    $$
    then \eqref{wavelet} is a Banach frame for $\dot{F}^s_{p,q}(\R)$ in the sense
    of \eqref{waveletexp}\,.
\end{description}
\end{Theorem}

\bproof Let us prove (a). First of all, we apply Theorem
\ref{coohom}/(i). Afterwards, we use Proposition \ref{rti} in order
to estimate the weight $w_{Y}(x,t)$ for
$Y=\dotspace{L}{s+d/2-d/q}$\,. We obtain 
\be\nonumber
  \begin{split}
    w_Y(x,t) &= \max\{t^{d(1/p-1/2)-s},t^{s-d(1/p-1/2)},t^{s+d/2},t^{-s+d/2}\}\\
    &\leq \left\{\begin{array}{rcl}
            t^{-r_1}&:&0<t<1\\
            t^{r_2}&:&t\geq 1.
    \end{array}\right.
  \end{split}
\ee
Let us distinguish the cases $s\geq 0$ and $s<0$. In the first case we can put
$r_1 = \max\{s-d(1/p-1/2),-s+d(1/p-1/2),s-d/2\}$, $r_2 = \max\{s+d/2,-s+d(1/p-1/2)\}$ and $v=0$.
Now we apply first Proposition \ref{propwiener}. This gives the condition
\be\label{s>0}
    0\leq s <  \min\{L,K\}-d(1-1/p)\,.
\ee
In the second case we put $r_1 = \max\{s-d(1/p-1/2),-s+d(1/p-1/2),-s-d/2\}$,
$r_2 = \max\{-s+d/2,-s+d(1/p-1/2)\}$ and $v=0$. With Proposition \ref{propwiener}
we obtain the condition
\be\label{s<0}
   -\min\{L,K\}+d/p<s<0\,.
\ee
Finally \eqref{s>0}, \eqref{s<0} and Theorem \ref{wbases} yield (a)\,.

{\em Step 2. } We prove (b). We apply Theorem \ref{coohom}/(iii) and
afterwards Proposition \ref{rti} and obtain for $Y =
\dotspace{P}{s+d/2-d/q,a}$ 
\be\nonumber
  \begin{split}
    w_Y(x,t) &= \max\{t^{d(1/p-1/2)-s},t^{s-d(1/p-1/2)},\\
    &~~~~~~~t^{s+d/2}\max\{1,t^{-a}\}(1+|x|^a),
    t^{-s+d/2}\max\{t^{-a},t^{a}\}(1+|x|)^a\}\\
    &\leq (1+|x|)^a\left\{\begin{array}{rcl}
            t^{-r_1}&:&0<t<1\\
            t^{r_2}&:&t\geq 1\,.
    \end{array}\right.
  \end{split}
\ee
First, we consider the case $s\geq 0$. We can put
$r_1 = \max\{s+a-d/2,s+d/2-d/p\}$, $r_2 = \max\{s+d/2,-s+d/2+a\}$ and $v=a$\,.
Proposition \ref{propwiener} gives the condition
$$
    0\leq s < \min\{L,K\}-\max\{a,d(1-1/p)\}
$$
which can be rewritten to
\be\nonumber
     0\leq s < \min\{L,K\}-d\max\{1/p,1/q,1-1/p\}
\ee
since $a$ can be chosen arbitrarily greater than $d\max\{1/p,1/q\}$\,. This gives the upper bound
in (b). Now we consider $s<0$. We put $r_1 = \max\{-s+d/p-d/2, s+d/2-d/p\}$,
$r_2 = -s+d/2+a$\,. This yields
$$
    -\min\{L,K\}+2a <s < 0
$$
and can be rewritten to
\be\nonumber
    -\min\{L,K\}+2d\max\{1/p,1/q\} < s <0\,.
\ee
This yields the lower bound in (b) and we are done. \eproof\\
\noindent
The following corollary is a consequence of Theorem \ref{wavdec} and the facts in Section \ref{appwavelets}.

\begin{corollary} Let $m>0$ and $(\Psi^0,\Psi^1) = (\varphi_m,\psi_m)$ the spline wavelet system of order $m$. Then
\begin{itemize}
    \item[(a)] If $1\leq p,q \leq \infty$ and
    $$
            -m+1+\frac{d}{p}<s<m-1-d\Big(1-\frac{1}{p}\Big)
    $$
    then \eqref{wavelet} is a Banach frame for $\dot{B}^s_{p,q}(\R)$ in the sense
    of \eqref{waveletexp}.
    \item[(b)] If $1\leq p <\infty$, $1 \leq q\leq \infty$, and
    $$
            -m+1+2d\max\Big\{\frac{1}{p},\frac{1}{q}\Big\}<s<m-1-d\max\Big\{\frac{1}{p},\frac{1}{q},1-\frac{1}{p}\Big\}
    $$
    then \eqref{wavelet} is a Banach frame for $\dot{F}^s_{p,q}(\R)$ in the sense
    of \eqref{waveletexp}\,.
\end{itemize}
\end{corollary}

\begin{remark}
 The (optimal) smoothness conditions in \cite{Bo95} are slightly weaker than (a) in case $d=1$. However, 
 compared to the approach of Triebel \cite{Tr06,Tr08}, we admit
 some more degree of freedom. The wavelet or atom does not have to be compactly supported. Additionally, in case $d=1$ we do
 not need that $\psi \in C^u(\re)$ where $u>s$. Indeed, the conditions in (a) and (b) are slightly weaker.
\end{remark}
\begin{remark}
  More examples can be obtained by using compactly supported Daubechies wavelets of a certain order or Meyer wavelets. Based on the
  underlying abstract result in Theorem \ref{wbases} even biorthogonal
  wavelet systems providing sufficiently high smoothness and vanishing moments are suitable for this issue. 
\end{remark}

\setcounter{section}{0}
\renewcommand{\thesection}{\Alph{section}}
\renewcommand{\theequation}{\Alph{section}.\arabic{equation}}
\renewcommand{\theTheorem}{\Alph{section}.\arabic{Theorem}}
\renewcommand{\thesubsection}{\Alph{section}.\arabic{subsection}}

\section{Appendix: Wavelets}

\subsection{The continuous wavelet transform}
\label{SectCWT}
The vector $g$ is said to be the analyzing vector for a function $f\in L_2(\R)$. The
continuous wavelet transform $W_gf$ is then defined
through
\begin{equation}\nonumber
    W_g f(x,t) = \langle T_x\mathcal{D}^{L_2}_t g, f\rangle\quad,\quad x\in \R, t>0\,,
\end{equation} where the bracket $\langle \cdot, \cdot \rangle$ denotes the
inner product in $L_2(\R)$. We can write it in terms of the convolution \eqref{conv} via 
\begin{equation}\label{CWTconv}
    W_g f(x,t) = [(\mathcal{D}^{L_2}_t g(-\cdot)) \ast \bar{f}](x) = t^{d/2}[(\mathcal{D}_t g(-\cdot)) \ast \bar{f}](x)\,.
\end{equation}
We call $g$ an admissible wavelet if
$$
    c_g:= \int\limits_{\R} \frac{|\cf g(\xi)|^2}{|\xi|^d}\,d\xi < \infty\,.
$$
If this is the case, then the family $\{T_x\mathcal{D}^{L_2}_t g\}_{t>0, x\in \R}$ represents a tight
continuous frame in $L_2(\re^d)$. For a proof
we refer to Theorem 1.5.1 in \cite{ra05-6}\,.\\

Let us now specify the conditions $(M_L)$, $(D)$, and $(S_K)$ which we intend to impose on functions $\Phi,\Psi\in L_2(\R)$ in order
to obtain a proper decay of the continuous wavelet transform $|W_{\Psi}\Phi(x,t)|$.

\begin{definition}\label{basedef} Let $L+1\in \n$, $K>0$ and fix the conditions $(D)$, $(M_L)$ and $(S_K)$ for
a function $\Psi \in L_2(\R)$.
\begin{enumerate}
     \item[$(D)$] For every $N\in \N$ there exists a constant
    $c_N$ such that
    $$
        |\Psi(x)| \leq \frac{c_N}{(1+|x|)^N}\,.
    $$
    \item[$(M_L)$] We have vanishing moments
    $$
        D^{\bar{\alpha}}\cf \Psi (0) = 0
    $$
    for all $|\bar{\alpha}|_1 \leq L$\,.
    \item[$(S_K)$] The function
    $$
        (1+|\xi|)^{K}|D^{\bar{\alpha}}\cf \Psi(\xi)|
    $$
    belongs to $L_1(\R)$ for every multi-index $\bar{\alpha}\in \n^d$.
\end{enumerate}
\end{definition}
\begin{remark} If a function $g \in L_2(\R)$ satisfies $(S_K)$ for some $K>0$
then by well-known properties of the Fourier transform we have $g
\in C^{\lfloor K \rfloor}(\R)$.
\end{remark}

The following lemma provides a useful decay result for the
continuous wavelet transform under certain smoothness, decay and
moment conditions, see also \cite{FrJa90, Ry99a, HeNe07} for similar results 
in a different language. It represents a continuation of \cite[Lem.
1]{Ry99a} where one deals with $\mathcal{S}(\R)$-functions
\begin{lemma}\label{help1} Let $L\in \n$, $K>0$ and $\Phi,\Psi, \Phi_0 \in L_2(\R)$.
   \begin{description}
   \item(i) Let $\Phi$ satisfy $(D)$, $(M_{L-1})$ and let $\Phi_0$ satisfy $(D)$, $(S_{K})$. Then for every
   $N\in \N$ there exists a constant $C_N$ such that the
   estimate \begin{equation}\label{f18}
       |(W_\Phi \Phi_0)(x,t)| \leq C_N \frac{t^{\min\{L,K\}+d/2}}{(1+|x|)^N}
   \end{equation}
   holds true for $x\in \R$ and $0<t<1$\,.
   \item(ii) Let $\Phi,\Psi$ satisfy $(D)$, $(M_{L-1})$ and $(S_K)$. For every $N\in \N$ there exists a constant $C_N$ such that the
   estimate 
   \begin{equation}\nonumber
       |(W_\Phi \Psi)(x,t)| \leq C_N
       \frac{t^{\min\{L,K\}+d/2}}{(1+t)^{2\min\{L,K\}+d}}\Big(1+\frac{|x|}{1+t}\Big)^{-N}
   \end{equation}
   holds true for $x\in \R$ and $0<t<\infty$\,.
   \end{description}
\end{lemma}

\bproof {\em Step 1.} Let us prove (i).
We follow the proof of Lemma 1 in \cite{Ry99a}. This reference deals
with $\mathcal{S}(\R)$-functions, which makes the situation much more easy. Assume without loss of generality 
$\Phi = \Phi(-\cdot)$ and $\Phi_0$ to be real-valued. Formula \eqref{CWTconv} gives
\begin{equation}\label{eq-21}
    |(W_{\Phi}\Phi_{0})(x,t)| = t^{d/2}|[(\mathcal{D}_t \Phi)\ast
    \Phi_0](x)|\,.
\end{equation}
Fix $0<t<1$. Obviously, the convolution
$(\mathcal{D}_t \Phi)\ast \Phi_0$ satisfies $(D)$. By well-known
properties of the Fourier transform the derivative
$D^{\bar{\alpha}}\cf((\mathcal{D}_t \Phi)\ast \Phi_0))(\xi)$ exists
for every multi-index $\bar{\alpha} \in \n^d$. For fixed
$\bar{\alpha}$ we estimate by using Leibniz' formula
\begin{equation}\label{eq1}
  \begin{split}
    |D^{\bar{\alpha}}\cf((\mathcal{D}_t \Phi)\ast \Phi_0))(\xi)|
    &= |D^{\bar{\alpha}}(\cf\Phi(t\xi)\cdot \cf \Phi_0 (\xi))|\\
    &\leq c_{\bar{\alpha}}\sum\limits_{\bar{\beta}\leq \bar{\alpha}}
    t^{|\bar{\beta}|_1}|D^{\bar{\beta}}\cf \Phi(t\xi)\cdot
    D^{\bar{\alpha}-\bar{\beta}}\cf \Phi_0(\xi)|\\
    &\leq c_{\bar{\alpha}}'t^L(1+|\xi|)^L\sum\limits_{\bar{\beta}\leq \bar{\alpha}}
    |D^{\bar{\alpha}-\bar{\beta}}\cf \Phi_0(\xi)|\,.
  \end{split}
\end{equation}
In the last step we used property $(M_{L-1})$. Assuming $K\geq L$ and exploiting $(S_K)$ we obtain that the left-hand
side of \eqref{eq1} belongs to $L_1(\R)$ and
\begin{equation}\label{eq-20}
    \|D^{\bar{\alpha}}\cf((\mathcal{D}_t \Phi)\ast
    \Phi_0))(\xi)|L_1(\R)\| \leq  c''_{\bar{\alpha}}\,t^{L}\,.
\end{equation}
We proceed as follows
\begin{equation}\nonumber
    \begin{split}
        \max\limits_{|\bar{\alpha}|_1\leq N+1}\|D^{\bar{\alpha}}\cf((\mathcal{D}_t \Phi)\ast
        \Phi_0))(\xi)|L_1(\R)\| &\geq
        \max\limits_{|\bar{\alpha}|_1\leq N+1}
        \|\cf^{-1}[D^{\bar{\alpha}}\cf((\mathcal{D}_t \Phi)\ast
        \Phi_0)]|L_{\infty}(\R)\|\\
        &\geq c_N\|(1+|x|)^N [(\mathcal{D}_t \Phi)\ast \Phi_0](x)|L_{\infty}(\R)\|\,.
    \end{split}
\end{equation}
This estimate together with \eqref{eq-21} and \eqref{eq-20} yields \eqref{f18}.\\
Let us finally assume
$K<L$ and return to \eqref{eq1}. Clearly, the resulting
inequality remains valid if we replace the exponent $L$
by $L'\in \n$ with $L' \leq L$. It is even possible to extend \eqref{eq1} to every
$0\leq L'' < L$ by the following argument. Let $L'' \notin \N$. We have on the one hand
$$
   \mbox{LHS}\eqref{eq1} \leq c'_{\bar{\alpha}}t^{\lfloor L'' \rfloor}(1+|\xi|)^{\lfloor L'' \rfloor}G(\xi)
$$
and on the other hand
$$
   \mbox{LHS}\eqref{eq1} \leq c'_{\bar{\alpha}}t^{\lfloor L'' \rfloor+1}(1+|\xi|)^{\lfloor L'' \rfloor+1}G(\xi)\,,
$$
where $G(\xi) = \sum_{\bar{\beta}\leq \bar{\alpha}}|D^{\bar{\alpha}-\bar{\beta}}\cf \Phi_0(\xi)|$\,. Choosing $0<\theta <1$
such that $L'' = (1-\theta)\lfloor L'' \rfloor + \theta(\lfloor L'' \rfloor+1)$ we obtain by a kind of interpolation argument
\begin{equation}\nonumber
  \begin{split}
  \mbox{LHS}\eqref{eq1} &= \mbox{LHS}\eqref{eq1}^{1-\theta}\mbox{LHS}\eqref{eq1}^{\theta}\\
  &\leq c'_{\bar{\alpha}} t^{L''}(1+|\xi|)^{L''}G(\xi)\,.
  \end{split}
\end{equation}
In particular, we obtain instead of $\eqref{eq1}$
$$
   |D^{\bar{\alpha}}\cf((\mathcal{D}_t \Phi)\ast \Phi_0))(\xi)| \leq
   c'_{\bar{\alpha}}t^{K}(1+|\xi|)^{K}\sum\limits_{\bar{\beta}\leq \bar{\alpha}}
    |D^{\bar{\alpha}-\bar{\beta}}\cf \Phi_0(\xi)|\quad,\quad \xi\in \R\,.
$$
We exploit property $(S_K)$ for $\Phi_0$ and proceed analogously as above. This proves \eqref{f18}.\frei\\\newline
{\em Step 2.} The estimate in (ii) is an immediate consequence of \eqref{f18} and the fact
$$
      (W_{\Phi}\Psi)(x,t) = (W_{\Psi}\Phi)(-x/t,1/t)\,.
$$
This completes the proof.\eproof

\begin{corollary}\label{remhelp1} Let $\Phi,\Psi$ belong to the Schwartz space $\mathcal{S}_0(\R)$. By Lemma \ref{help1}/(ii) for every $L,N\in \N$ there is a constant
$C_{L,N}>0$ such that
$$
    |(W_{\Phi}\Psi)(x,t)| \leq C_{L,N}\frac{t^{L+d/2}}{(1+t)^{2L+d}}\Big(1+\frac{|x|}{1+t}\Big)^{-N}\quad,\quad x\in \R, t>0\,.
$$
Additionally, we obtain for $\Phi\in \mathcal{S}_0(\R)$ and $\Phi_0\in \mathcal{S}(\R)$ that
$$
    |(W_{\Phi}\Phi_0)(x,t)| \leq C_{L,N}\frac{t^{L+d/2}}{(1+|x|)^N}\quad,\quad x\in \R, 0<t<1\,.
$$

\end{corollary}

\subsection{Orthonormal wavelet bases}\label{appwavelets}
The following
Lemma is proved in Wojtaszczyk \cite[5.1]{Wo97}.
\begin{lemma}\label{dwavelet} Suppose we have a multiresolution analysis
in $L_2(\re)$ with scaling functions $\Psi^{0}(t)$ and associated
wavelets $\Psi^{1}(t)$. Let $E = \{0,1\}^d\setminus
\{(0,...,0)\}$. For $c = (c_1,...,c_d) \in E$ let $\Psi^c =
\bigotimes_{j=1}^d \Psi^{c_j}$. Then the system
\be\label{wavelet}
        \Big\{2^{\frac{jd}{2}}\Psi^{c}(2^jx-k)\Big\}_{c\in E,j\in \zz,k\in \Z}
\ee is an orthonormal basis in $L_2(\R)$.
\end{lemma}

\subsubsection*{Spline wavelets}
As a main example we will consider the spline wavelet system.  The normalized cardinal B-spline of order
$m+1$ is given by
\[
\cn_{m+1} (x):= \cn_m * \cx (x)\, , \qquad x \in \re\, , \quad m \in
\N\, ,
\]
beginning with $\cn_1 = \cx$, the characteristic function of the interval $(0,1)$.
By
\[
\varphi_m (x):= \frac{1}{\sqrt{2\pi}} \, \cfi \Big[\frac{\cf \cn_m
(\xi)}{\Big(\sum\limits_{k=-\infty}^\infty |\cf \cn_m (\xi + 2\pi
k)|^2\Big)^{1/2}}\Big](x)\, , \qquad x \in \re\, ,
\]
we obtain an orthonormal scaling function which is again a spline of
order $m$.
Finally, by
\[
\psi_m (x) := \sum_{k=-\infty}^\infty \langle \, \varphi_m (t/2),
\varphi_m (t-k)\rangle\, (-1)^k \, \varphi_m (2x+k+1)
\]
the generator of an orthonormal wavelet system is defined. For $m=1$
it is easily checked that $-\psi_1 (x-1)$ is the Haar wavelet. In
general these functions $\psi_m$ have the following properties:
\begin{itemize}
\item $\psi_m$ restricted to intervals $[\frac{k}{2},\frac{k+1}{2}]$, $k\in \zz$,
is a polynomial of degree at most $m-1$.
\item $\psi_m \in C^{m-2} (\re)$ if $m\ge 2$.
\item $\psi_m^{(m-2)}$ is uniformly Lipschitz continuous on
$\re$ if $m \ge 2$.
\item
The function $\psi_m$ satisfies a moment condition of order $m-1$,
i.e.
\[
\int_{-\infty}^\infty x^\ell \, \psi_m (x)\, dx =0\, , \qquad \ell
=0,1,\ldots, m-1\, .
\]
In particular, $\psi_m$ satisfies $(M_{L})$ for $0<L\leq m$ and $\psi_m, \varphi_m$ satisfy $(D)$ and $(S_{K})$ for $K<m-1$.
\end{itemize}


\begin{thebibliography}{10}

\bibitem{Ba03}
D.~B. {B}azarkhanov.
\newblock {C}haracterizations of the {N}ikol'skij-{B}esov and
  {L}izorkin-{T}riebel function spaces of mixed smoothness.
\newblock {\em Proc. Steklov Inst. Math.}, 243:53--65, 2003.

\bibitem{Bo95}
G.~{B}ourdaud.
\newblock Ondelettes et espaces de {B}esov.
\newblock {\em Rev. Matem. Iberoam.}, 11:477--512, 1995.

\bibitem{BuPaTa96}
H.-Q. {B}ui, M.~{P}aluszy\'nski, and M.~H. {T}aibleson.
\newblock A maximal function characterization of weighted {B}esov-{L}ipschitz
  and {T}riebel-{L}izorkin spaces.
\newblock {\em Stud. Math.}, 119(3):219--246, 1996.

\bibitem{BuPaTa97}
H.-Q. {B}ui, M.~{P}aluszy\'nski, and M.~H. {T}aibleson.
\newblock Characterization of the {B}esov-{L}ipschitz and {T}riebel-{L}izorkin
  spaces, the case $q<1$.
\newblock {\em J. Four. Anal. and Appl. (special issue)}, 3:837--846, 1997.

\bibitem{CoMeSt85}
R.~R. {C}oifman, Y.~{M}eyer, and E.~M. {S}tein.
\newblock {S}ome new function spaces and their application to harmonic
  analysis.
\newblock {\em Journ. of Funct. Anal.}, 62:304--335, 1985.

\bibitem{FeSt71}
C.~{F}efferman and E.~M. {S}tein.
\newblock {S}ome maximal inequalities.
\newblock {\em Amer. Journ. Math.}, 93:107--115, 1971.

\bibitem{FeGr86}
H.~G. Feichtinger and K.~Gr{\"o}chenig.
\newblock A unified approach to atomic decompositions via integrable group
  representations.
\newblock In {\em Function spaces and applications (Lund, 1986)}, volume 1302
  of {\em Lecture Notes in Math.}, pages 52--73. Springer, Berlin, 1988.

\bibitem{FeGr89a}
H.~G. {F}eichtinger and K.~{G}r{\"o}chenig.
\newblock {B}anach spaces related to integrable group representations and their
  atomic decompositions, {I}.
\newblock {\em Journ. Funct. Anal.}, 21:307--340, 1989.

\bibitem{FeGr89b}
H.~G. {F}eichtinger and K.~{G}r{\"o}chenig.
\newblock {B}anach spaces related to integrable group representations and their
  atomic decompositions, {II}.
\newblock {\em {M}onatsh. {M}athem.}, 108:129--148, 1989.

\bibitem{fegr92-1}
H.~G. {F}eichtinger and K.~{G}r{\"o}chenig.
\newblock {G}abor wavelets and the {H}eisenberg group: {G}abor expansions and
  short time {F}ourier transform from the group theoretical point of view.
\newblock In C.~K. {C}hui, editor, {\em {W}avelets :a tutorial in theory and
  applications}, volume~2 of {\em {W}avelet {A}nal. {A}ppl.}, pages 359--397.
  {A}cademic {P}ress, 1992.

\bibitem{fora05}
M.~{F}ornasier and H.~Rauhut.
\newblock {C}ontinuous frames, function spaces, and the discretization problem.
\newblock {\em {J}. {F}ourier {A}nal. {A}ppl.}, 11(3):245--287, 2005.

\bibitem{FrJa90}
M.~{F}razier and B.~{J}awerth.
\newblock {A} discrete transform and decompositions of distribution spaces.
\newblock {\em Journ. of Funct. Anal.}, 93:34--170, 1990.

\bibitem{Gr88}
K.~Gr{\"o}chenig.
\newblock Unconditional bases in translation and dilation invariant function
  spaces on {$\mathbb{R}\sp n$}.
\newblock In {\em Constructive theory of functions (Varna, 1987)}, pages
  174--183. Publ. House Bulgar. Acad. Sci., Sofia, 1988.

\bibitem{Gr91}
K.~{G}r{\"o}chenig.
\newblock {D}escribing functions: atomic decompositions versus frames.
\newblock {\em {M}onatsh. {M}athem.}, 112:1--41, 1991.

\bibitem{HaDiss10}
M.~{H}ansen.
\newblock {\em {N}onlinear {A}pproximation and {F}unction {S}paces of
  {D}ominating {M}ixed {S}moothness}.
\newblock PhD thesis, FSU Jena, Germany, 2010.

\bibitem{HeNe07}
L.~I. {H}edberg and Y.~{N}etrusov.
\newblock {A}n axiomatic approach to function spaces, spectral synthesis, and
  {L}uzin approximation.
\newblock {\em Mem. Amer. Math. Soc.}, 188(882), 2007.

\bibitem{Ke09}
H.~{K}empka.
\newblock 2-microlocal spaces of variable integrability.
\newblock {\em {R}ev. {M}at. {C}omplut.}, 22(1):227--251, 2009.

\bibitem{Ov84}
V.~I. {O}vchinnikov.
\newblock {\em {T}he {M}ethod of {O}rbits in {I}nterpolation {T}heory},
  volume~1 of {\em {M}athematical {R}eports}.
\newblock {H}arwood {A}cad. {P}ublishers, 1984.

\bibitem{Pe75}
J.~{P}eetre.
\newblock {O}n spaces of {T}riebel-{L}izorkin type.
\newblock {\em {A}rk. {M}at.}, 13:123--130, 1975.

\bibitem{ra05-6}
H.~{R}auhut.
\newblock {\em {T}ime-frequency and wavelet analysis of functions with symmetry
  properties}.
\newblock {L}ogos-{V}erlag, 2005.
\newblock {P}h{D} thesis.

\bibitem{ra05-3}
H.~{R}auhut.
\newblock {C}oorbit space theory for quasi-{B}anach spaces.
\newblock {\em {S}tudia {M}ath.}, 180(3):237--253, 2007.

\bibitem{ra05-4}
H.~{R}auhut.
\newblock {W}iener amalgam spaces with respect to quasi-{B}anach spaces.
\newblock {\em {C}olloq. {M}ath.}, 109(2):345--362, 2007.

\bibitem{Ry99a}
V.~S. {R}ychkov.
\newblock {O}n a theorem of {B}ui, {P}aluszy\'nski and {T}aibleson.
\newblock {\em {P}roc. {S}teklov {I}nst.}, 227:280--292, 1999.

\bibitem{Ry99b}
V.~S. {R}ychkov.
\newblock {O}n restrictions and extensions of the {B}esov and
  {T}riebel-{L}izorkin spaces with respect to {L}ipschitz domains.
\newblock {\em Journ. Lond. Math. Soc.}, 60:237--257, 1999.

\bibitem{Ry01}
V.~S. {R}ychkov.
\newblock {L}ittlewood-{P}aley theory and function spaces with
  ${A}_p^{\mbox{loc}}$ weights.
\newblock {\em {M}ath. {N}achr.}, 224:145--180, 2001.

\bibitem{StWe71}
E.~M. Stein and G.~{W}eiss.
\newblock {\em {I}ntroduction to {F}ourier {A}nalysis on {E}uclidean {S}paces}.
\newblock Princeton Univ. Press, 1971.

\bibitem{StTo89}
J.-O. Str{\"o}mberg and A.~{T}orchinsky.
\newblock {\em {W}eighted {H}ardy {S}paces}.
\newblock {L}ecture {N}otes in {M}athematics, {V}ol. 1381. {S}pringer, 1989.

\bibitem{Tr83}
H.~{T}riebel.
\newblock {\em {T}heory of {F}unction {S}paces}.
\newblock {B}irkh{\"a}user, Basel, 1983.

\bibitem{Tr88}
H.~{T}riebel.
\newblock {C}haracterizations of {B}esov-{H}ardy-{S}obolev spaces: a unified
  approach.
\newblock {\em {J}ourn. of {A}pprox. {T}heory}, 52:162--203, 1988.

\bibitem{Tr92}
H.~{T}riebel.
\newblock {\em {T}heory of {F}unction {S}paces II}.
\newblock {B}irkh{\"a}user, Basel, 1992.

\bibitem{Tr06}
H.~{T}riebel.
\newblock {\em Theory of {F}unction {S}paces III}.
\newblock Birkh{\"a}user, Basel, 2006.

\bibitem{Tr08}
H.~{T}riebel.
\newblock {\em Function {S}paces and {W}avelets on {D}omains}.
\newblock EMS Publishing House, Z{\"u}rich, 2008.

\bibitem{Vy06}
J.~{V}ybiral.
\newblock {F}unction spaces with dominating mixed smoothness.
\newblock {\em Diss. Math.}, 436:1--73, 2006.

\bibitem{Wo97}
P.~{W}ojtaszczyk.
\newblock {\em {A} {M}athematical {I}ntroduction to {W}avelets}.
\newblock {C}ambridge {U}niversity {P}ress, 1997.

\end{thebibliography}
\end{document}